\RequirePackage{fix-cm}
\documentclass[smallextended]{svjour3}       
\smartqed  
\usepackage{amssymb}
\usepackage{graphicx}
\usepackage{amsmath}
\usepackage{multirow,url}
\usepackage{float}
 \usepackage{rotating}
 \usepackage{mathabx}
\usepackage{mathptmx}      
\usepackage{xcolor}
\usepackage{latexsym}
\usepackage{hyperref}
\hypersetup{hypertex=true,
colorlinks=true,
linkcolor=blue,
anchorcolor=blue,
citecolor=blue}
\newtheorem{assumption}{Assumption}

\begin{document}

\title{A preconditioned MINRES method for block lower triangular Toeplitz systems}
\author{\small Congcong Li\and Xuelei Lin \and Sean Hon \and Shu-Lin Wu
}


\institute{Congcong Li \at School of Science, Harbin Institute of Technology,	Shenzhen 518055, China 
              \at Department of Mathematics, Hong Kong Baptist University, Kowloon Tong, Hong Kong SAR \\             \email{22482245@life.hkbu.edu.hk}         
           \and
           Corresponding author: Xuelei Lin \at
              School of Science, Harbin Institute of Technology,	Shenzhen 518055, China\\
       \email{linxuelei@hit.edu.cn}
       \and
          Sean Hon \at
              Department of Mathematics, Hong Kong Baptist University, Kowloon Tong, Hong Kong SAR\\
        \email{seanyshon@hkbu.edu.hk}
        \and
        Shu-Lin Wu \at
              School of Mathematics and Statistics, Northeast Normal University, Changchun 130024, China\\
        \email{wushulin84@hotmail.com}
}

\date{Received: date / Accepted: date}

\maketitle

\begin{abstract}
In this study, a novel preconditioner based on the absolute-value block $\alpha$-circulant matrix approximation is developed, specifically designed for nonsymmetric dense block lower triangular Toeplitz (BLTT) systems that emerge from the numerical discretization of evolutionary equations. Our preconditioner is constructed by taking an absolute-value of a block $\alpha$-circulant matrix approximation to the BLTT matrix. To apply our preconditioner,  the original BLTT linear system is converted into a symmetric form by applying a time-reversing permutation transformation. Then, with our preconditioner, the preconditioned minimal residual method (MINRES) solver is employed to solve the symmetrized linear system. With properly chosen $\alpha$, the eigenvalues of the preconditioned matrix are proven to be clustered around $\pm1$ without any significant outliers. With  the clustered spectrum, we show  that the preconditioned MINRES solver for the preconditioned system has a convergence rate independent of system size. The efficacy of the proposed preconditioner is corroborated by our numerical experiments, which reveal that it attains optimal convergence.
\end{abstract}

\section{Introduction}
\label{intro}
Consider a  BLTT system, given by:
\begin{equation}\label{blttsystem}
  \mathcal{A} \mathbf{u}=\mathbf{f},  
\end{equation}
where $\mathbf{u}$ is the unknown vector, $\mathbf{f}$ is the right hand side, and
\begin{equation}\label{amatdef}
\mathcal{A}=\left[\begin{array}{ccccc}
A_{(0)} & & & & \\
A_{(1)} & \ddots & & & \\
\vdots & \ddots & \ddots & & \\
\vdots & \ddots & \ddots & \ddots & \\
A_{(N-1)} & \cdots & \cdots & A_{(1)} & A_{(0)}
\end{array}\right]
\end{equation}
is an $MN$ $\times$ $MN$ BLTT matrix with symmetric blocks $A_{(k)}\in\mathbb{R}^{M\times M}$, $k=0, \dots, N-1$. The BLTT system \eqref{blttsystem}  commonly referred  as the all-at-once system, is obtained through time-space discretization of evolutionary equations containing classical or fractional temporal derivatives. These include the initial-boundary value problem of heat equations, wave equations, evolutionary convection diffusion
equations as well as fractional sub-diffusion equations, as reported in \cite{liu2020fast,mcdonald2018preconditioning,lin2021all,lin2021parallel,ke2015fast,hon2023block,hon2022sine,gu2017fast,lu2015fast,lin2016fast,lin2018separable}.

In general, Krylov subspace methods are highly appropriate for solving linear systems with Toeplitz structure such as the one described in \eqref{blttsystem}. Furthermore, it is worth mentioning that incorporating a preconditioner into the iterative procedure can lead to faster convergence, fewer iterations, and improved robustness for ill-conditioned matrices \cite{chan1996conjugate}. This can make preconditioned Krylov subspace methods a more practical and efficient approach for solving large linear systems of equations. The use of methods such as GMRES \cite{saad1986gmres}, QMR \cite{freund1991qmr}, SQMR \cite{freund1995software}, and BiCG-STAB \cite{van1992bi} for solving systems of equations with nonsymmetric matrices is common. However, it is important to note that these methods do not simultaneously minimize a relevant quantity and may not have short-term recurrences for general matrices. Additionally, their convergence cannot be bounded using only the eigenvalues or singular values \cite{greenbaum1996any}, which is in contrast to symmetric Krvlov subspace methods designed for symmetric systems, such as the preconditioned conjugate gradient (PCG) method \cite{hestenes1952methods} and the preconditioned MINRES method \cite{paige1975solution}. Consequently, it is challenging to develop a widely applicable convergence theory for many Krylov iterative methods in nonsymmetric linear systems.

To tackle this issue, Pestana and Wathen in \cite{pestana2015precon} developed a reordering technique to convert a nonsymmetric scalar Toeplitz system into  an indefinite symmetric Hankel system and a preconditioned MINRES method was proposed to solve the symmetrized system. The fast convergence of the preconditioned MINRES method proposed in \cite{pestana2015precon} is based on orthogonal-plus-low-rank-plus-small-norm decomposition of the preconditioned matrix. To extend this preconditioned MINRES method to solving BLTT system, McDonald et al. \cite{mcdonald2018preconditioning} developed a time-reversing permutation transformation for symmetrizing BLTT system and proposed an absolute-value block circulant preconditioner to accelerate the convergence of the MINRES method for solving the symmetrized BLTT system arising from heat equations. However, such a preconditioning technique can not keep its efficiency in solving the BLTT system \eqref{blttsystem} as rank of the low-rank part would be enlarged by the block size $M$.  Lin and Hon \cite{lin2023block} proposed a preconditioned MINRES method based on the absolute value block $\alpha$-circulant matrix for a sparse BLTT system arising from solving wave equations and showed that the preconditioned MINRES solver has a convergence rate independent of matrix size. However, the theoretical results can not be applied to the dense BLTT linear system \eqref{blttsystem}. To conclude, the study of symmetric Krylov subspace methods for the BLTT system is still at its infancy. 

To fill this gap, we adopt the time-reversing permutation matrix proposed in \cite{McDonald2017,mcdonald2018preconditioning} to convert the dense BLTT linear system \eqref{blttsystem} into a symmetric indefinite system \eqref{symmetrizedblttsys}. Then, an absolute-value block $\alpha$-circulant (ABAC) preconditioner is proposed to accelerate the convergence of MINRES solver for the symmetrized system. Our ABAC preconditioner is actually a generalization of the absolute-value block circulant preconditioner proposed in \cite{mcdonald2018preconditioning} since the ABAC preconditioner with $\alpha=1$ is exactly the absolute-value block circulant preconditioner. Moreover, we show theoretically that the preconditioned MINRES with ABAC preconditioner has a convergence rate independent of $M$ and $N$ for the symmetrized dense BLTT linear system \eqref{symmetrizedblttsys} when the value of $\alpha$ is properly chosen.
Numerical results are reported to demonstrate the efficiency of the proposed ABAC preconditioner. 
 
 The rest of the paper is organized as follows. In section \ref{nblttsymmetrization}, we present the symmetrization of the BLTT system and introduce some assumptions for the blocks $A_{(k)}$ ($0\leq k\leq N-1$). In section \ref{abacpre}, we present the derivation of the ABAC preconditioner. In section \ref{convanaly}, the convergence rate of preconditioned MINRES method with the proposed ABAC preconditioner for the symmetrized linear system is analyzed. In section \ref{numexams},  we introduce the implementation details of the matrix-vector product associated with the preconditioned matrix and present the numerical results of  the proposed solver. At last, the concluding remarks are given in section \ref{conclusions}.

\section{Symmetrization of the BLTT system \eqref{blttsystem}}\label{nblttsymmetrization}
As mentioned in the introduction, instead of solving the BLTT system, we firstly apply a time-reversing permutation matrix to convert \eqref{blttsystem} into a symmetric system as follows
\begin{equation}\label{symmetrizedblttsys}
	\mathcal{Y} \mathcal{A} \mathbf{u}=\mathcal{Y} \mathbf{f},
\end{equation}
where
\begin{equation*}
	\mathcal{Y} \mathcal{A}=\left[\begin{array}{ccccc} 
		A_{(N-1)} & \cdots & \cdots & A_{(1)} & A_{(0)} \\
		\vdots& \reflectbox{$\ddots$} & \reflectbox{$\ddots$} & \reflectbox{$\ddots$} & \\
		\vdots & \reflectbox{$\ddots$} &\reflectbox{$\ddots$}  & & \\
		A_{(1)}  & \reflectbox{$\ddots$} & & & \\
		A_{(0)}  & & & &
	\end{array}\right],
\end{equation*}
$\mathcal{Y}=Y_N \otimes \mathbf{I}_M$ with $Y_N \in \mathbb{R}^{N \times N}$ being the anti-identity matrix, i.e., 
\begin{equation*}
Y_N=\left[\begin{array}[c]{ccc}
&&1\\
&\reflectbox{$\ddots$} &\\
1&&
\end{array}\right]\in\mathbb{R}^{N\times N}.
\end{equation*}
${\bf I}_{k}$ denotes the $k\times k$ identity throughout this paper. Clearly, the new coefficient matrix $\mathcal{Y} \mathcal{A}$ is symmetric since the blocks $A_{(k)}$ are real symmetric. $\mathcal{Y}$ is the so-called time-reversing permutation matrix. The eigenvalues of the symmetrized matrix $\mathcal{Y} \mathcal{A}$ have been shown in previous studies \cite{mazza2019spectral,hon2019note,ferrari2019eigenvalue} to be asymptotically distributed as $\pm|\hat{g}|$, where $\hat{g}$ represents the spectral symbol of $\mathcal{A}$. It can be inferred from such distribution that $\mathcal{Y} \mathcal{A}$ is generally an indefinite matrix.   

\begin{definition}\label{normalmatabsdef}
For a normal matrix ${\bf H}\in\mathbb{C}^{K\times K}$, its absolute value is defined as $|{\bf H}|:={\bf Q}^{\rm *}{\rm diag}(|\lambda_i|)_{i=1}^{K}{\bf Q}$, where ${\bf H}={\bf Q}^{\rm *}{\rm diag}(\lambda_i)_{i=1}^{K}{\bf Q}$ is the unitary diagonalization form of ${\bf H}$.
\end{definition}

 For a Hermitian matrix ${\bf H}$, denote by $\lambda_{\min}({\bf H})$ and $\lambda_{\max}({\bf H})$, the minimal and the maximal eigenvalues of ${\bf H}$, respectively.

In general, \eqref{symmetrizedblttsys} is indefinite. In the lateral section, we will propose the ABAC preconditioner and apply it as a preconditioner of preconditioned MINRES to solve \eqref{symmetrizedblttsys}. To prepare for the convergence analysis of the preconditioned MINRES solver, we list assumptions on the real symmetric blocks $A_{(k)}$ ($0\leq k\leq N-1$) as follows
\begin{assumption}\label{akassumps}
\item[(i)] The blocks $A_{(k)}$ commute with one another;
\item[(ii)]There exists a positive constant $c_0$ independent of both $N$ and $M$ such that
\begin{equation*}
	\lambda_{\min }\left(A_{(0)}-\sum_{i=1}^{N-1}\left|A_{(i)}\right|\right) \geq c_0 > 0.
\end{equation*}
\end{assumption}

\begin{remark}\label{assumptionakremark}
By combining Assumption \ref{akassumps}${\bf (i)}$ with the real symmetry of $A_{(k)}$'s,  it can be inferred that $A_{(k)}$'s are simultaneously diagonalizable by an orthogonal matrix ${\bf U}\in\mathbb{R}^{M\times M}$.
    Assumption \ref{akassumps}${\bf (i)}$-${\bf (ii)}$ are satisfied by many temporal discretization schemes combined with central difference/finite element discretization in space for evolutionary equations, such as backward Euler in time for heat equations and convection diffusion equations \cite{mcdonald2018preconditioning}, the finite difference discretization in time for fractional sub-diffusion equations \cite{lin2016fast,gao2011compact}, implicit leap-frog finite difference scheme in time for wave equations \cite{liu2020fast}. More examples of discretization of evolutionary equations satisfying Assumption \ref{akassumps} will be introduced in Section \ref{numexams}.
\end{remark}

\section{The derivation of ABAC preconditioner}\label{abacpre}
In this section, we will present the derivation of the ABAC preconditioner.

Denote 
\begin{equation*}
H_\alpha^{(k)}=\left[\begin{array}{cc}0 & \alpha \mathbf{I}_k \\ \mathbf{I}_{N-k} & 0\end{array}\right] \in \mathbb{R}^{N \times N},~ 1\leq k\leq N-1.
\end{equation*}
Then, it is straightforward to verify that
\begin{equation}\label{eq:tensorsumblttdef2}
\mathcal{A}  =\mathbf{I}_N \otimes A_{(0)}+\sum_{k=1}^{N-1} H_0^{(k)} \otimes A_{(k)}. 
\end{equation}
In \cite{mcdonald2018preconditioning}, an absolute-value block circulant preconditioner is proposed for the symmetrized BLTT system \eqref{symmetrizedblttsys}, which is defined by $|{\bf C}|$ with the normal matrix ${\bf C}$ given by
\begin{equation*}
{\bf C}=\mathbf{I}_N \otimes A_{(0)}+\sum_{k=1}^{N-1} H_1^{(k)} \otimes A_{(k)}.
\end{equation*}
Preconditioned MINRES with $|{\bf C}|$ as preconditioner for solving \eqref{symmetrizedblttsys} has a convergence rate dependent on $M$ (the larger $M$ is, the slower the convergence is); see, e.g., \cite{mcdonald2018preconditioning}. The drawback of the absolute-value block circulant preconditioner comes from the fact that $H_1^{(k)} $ is not a good approximation to $H_0^{(k)}$. To remedy this, we firstly consider the following approximation matrix
\begin{equation}\label{tensordef2calpha}
\mathcal{C}_{\alpha}:=\mathbf{I}_N \otimes A_{(0)}+\sum_{k=1}^{N-1} H_{\alpha}^{(k)} \otimes A_{(k)},\quad \alpha\in(0,1].
\end{equation}
Clearly, for sufficiently small $\alpha$, $H_{\alpha}^{(k)}$ is close to $H_0^{(k)} $ and thus $\mathcal{C}_{\alpha}$ is close to $\mathcal{A}$. Although $\mathcal{C}_{\alpha}$ is no longer normal for $\alpha\in(0,1)$, we can extend the concept of `absolute value' slightly in order to define an absolute value for $\mathcal{C}_{\alpha}$. 

\begin{definition}\label{diagonalizablematsqrtdef}
For a diagonalizable matrix ${\bf B}\in\mathbb{C}^{K\times K}$, let ${\bf B}={\bf X}^{-1}{\rm diag}(\lambda_i)_{i=1}^{K}{\bf X}$ be the diagonalization form of ${\bf B}$. Then, a square root of ${\bf B}$ is defined as $${\bf B}^{\frac{1}{2}}:={\bf X}^{-1}{\rm diag}(\lambda_i^{\frac{1}{2}})_{i=1}^{K}{\bf X},$$ where $\lambda_i^{\frac{1}{2}}$ is defined as the principle branch of complex square root of $\lambda_i$.
\end{definition}

\begin{definition}\label{diaglizablematabsdef}
For a diagonalizable matrix ${\bf B}\in\mathbb{C}^{K\times K}$, its absolute value is defined as
\begin{equation*}
|{\bf B}|:=(B^{\frac{1}{2}})^{\rm *}B^{\frac{1}{2}}.
\end{equation*}
\end{definition}
Clearly, Definition \ref{diaglizablematabsdef} is consistent with Definition \ref{normalmatabsdef}, i.e., for a normal matrix ${\bf H}$, its absolute values by the two definitions are the same. With Definition \ref{diaglizablematabsdef}, one can define the absolute value of any diagonalizable matrix. Fortunately, $\mathcal{C}_{\alpha}$ is diagonalizable. To see this, we firstly observe from Remark \ref{assumptionakremark} that there exists an orthonormal matrix $ U\in\mathbb{R}^{M\times M}$ such that 
\begin{equation*}
A_{(k)}=U\Lambda_k U^{\rm T},\quad \Lambda_k={\rm diag}(\lambda_i^{(k)})_{i=1}^{M}\in\mathbb{R}^{M\times M},\quad k=0,1,...,N-1.
\end{equation*}
With the above diagonalization formulas, one can rewrite $\mathcal{C}_{\alpha}$ as
\begin{equation*}
\mathcal{C}_{\alpha}=({\bf I}_N\otimes U)\left(\mathbf{I}_N\otimes \Lambda_0+\sum\limits_{k=1}^{N-1}H_{\alpha}^{(k)}\otimes \Lambda_k \right)({\bf I}_N\otimes U^{\rm T}).
\end{equation*}

Let $\Pi\in\mathbb{R}^{MN\times MN}$ be the following permutation matrix exchanging  the ordering of Kronecker product:
\begin{equation*}
\Pi({\bf A}\otimes {\bf B})\Pi^{\rm T}={\bf B}\otimes {\bf A},\quad \forall {\bf A}\in\mathbb{C}^{M\times M},~\forall{\bf B}\in\mathbb{C}^{N\times N}.
\end{equation*}

Then, it is clear that
\begin{equation}\label{calphapipermuteform}
\mathcal{C}_{\alpha}=({\bf I}_N\otimes U)\Pi\left(\Lambda_0\otimes\mathbf{I}_N+\sum\limits_{k=1}^{N-1}\Lambda_k\otimes  H_{\alpha}^{(k)} \right)\Pi^{\rm T}({\bf I}_N\otimes U^{\rm T}).
\end{equation}
It is straightforward to verify that
 \begin{align}
 &\Lambda_0\otimes\mathbf{I}_N+\sum\limits_{k=1}^{N-1}\Lambda_k\otimes  H_{\alpha}^{(k)}={\rm blockdiag}(\widehat{C}_{\alpha,i})_{i=1}^{M},\label{preblockdiagform}\\
 &\widehat{C}_{\alpha, i}=\left[\begin{array}{ccccc}
 	\lambda_i^{(0)} & \alpha \lambda_i^{(N-1)} & \cdots & \cdots & \alpha \lambda_i^{(1)} \\
 	\lambda_i^{(1)} & \ddots & \ddots & \ddots & \vdots \\
 	\vdots & \ddots & \ddots & \ddots & \vdots \\
 	\vdots & \ddots & \ddots & \ddots & \alpha \lambda_i^{(N-1)} \\
 	\lambda_i^{(N-1)} & \cdots & \cdots & \lambda_i^{(1)} & \lambda_i^{(0)}
 \end{array}\right] \in \mathbb{R}^{N \times N}.\label{hatcalphaidef}
 \end{align}
Observing from the structure of $\widehat{C}_{\alpha, i}$, we see that $\widehat{C}_{\alpha, i}$ is an $\alpha$-circulant matrix for each $1\leq i\leq M$. According to \cite{bertaccini2003block,bini2005numerical}, $\widehat{C}_{\alpha, i}$'s ($i=1,2,...,M$) are simultaneously diagonalizable as follows
\begin{equation}\label{hatcalphaidiagform}
	\widehat{C}_{\alpha, i}=D_\alpha^{-1} \mathbb{F} \Lambda_{\alpha, i} \mathbb{F}^* D_\alpha, \quad \Lambda_{\alpha, i}=\operatorname{diag}\left(\lambda_i^{(k, \alpha)}\right)_{k=1}^{N},
\end{equation}
where
\begin{align}
	D_\alpha&=\operatorname{diag}\left(\alpha^{\frac{i-1}{N}}\right)_{i=1}^N, \mathbb{F}=\frac{1}{\sqrt{N}}\left[\theta_N^{(i-1)(j-1)}\right]_{i, j=1}^N, \theta_N=\exp \left(\frac{2 \pi \mathbf{i}}{N}\right), \quad \mathbf{i}=\sqrt{-1}, \notag\\
	\lambda_i^{(k, \alpha)}&=\sum_{j=1}^N \widehat{C}_{\alpha, i}(j, 1) \alpha^{\frac{j-1}{N}} \theta_N^{-(k-1)(j-1)}=\left(\sqrt{N} \mathbb{F}^*D_{\alpha} \widehat{C}_{\alpha, i}(:, 1)\right)(k) .\label{hatcalphaieigvaldef}
\end{align}

Combining \eqref{calphapipermuteform}, \eqref{preblockdiagform} and \eqref{hatcalphaidiagform}, we see that $\mathcal{C}_{\alpha}$ is diagonalizable as follows
\begin{align}\label{calphadiagform}
\mathcal{C}_{\alpha}&=\left(\mathbf{I}_N \otimes U\right) \Pi \operatorname{blockdiag}\left(D_\alpha^{-1} \mathbb{F} \Lambda_{\alpha, i} \mathbb{F}^* D_\alpha\right)_{i=1}^M \Pi^{\rm T}\left(\mathbf{I}_N \otimes U\right)^{\rm T}\notag\\
&=[\left(\mathbf{I}_N \otimes U\right) \Pi [{\bf I}_{M}\otimes (D_{\alpha}^{-1}\mathbb{F})]]{\rm diag}(\lambda_i^{(k, \alpha)})_{k=1,i=1}^{N,M}[\left(\mathbf{I}_N \otimes U\right) \Pi [{\bf I}_{M}\otimes (D_{\alpha}^{-1}\mathbb{F})]]^{-1}.
\end{align}
With \eqref{calphadiagform}, Definition \ref{diaglizablematabsdef} gives an absolute-value of $\mathcal{C}_{\alpha}$ as follows
\begin{equation}\label{calphaabsdef}
\left|\mathcal{C}_\alpha\right|:=\left(\mathcal{C}_\alpha^{\frac{1}{2}}\right)^* \mathcal{C}_\alpha^{\frac{1}{2}},
\end{equation}
where
\begin{align}
\mathcal{C}_\alpha^{\frac{1}{2}}&=\left(\mathbf{I}_N \otimes U\right) \Pi \operatorname{blockdiag}\left(D_\alpha^{-1} \mathbb{F} \Lambda_{\alpha, i}^{\frac{1}{2}} \mathbb{F}^* D_\alpha\right)_{i=1}^M \Pi^{\rm T}\left(\mathbf{I}_N \otimes U\right)^{\rm T}\label{sqrtcalphabdiagform}\\
&=[\left(\mathbf{I}_N \otimes U\right) \Pi [{\bf I}_{M}\otimes (D_{\alpha}^{-1}\mathbb{F})]]{\rm diag}((\lambda_i^{(k, \alpha)})^{\frac{1}{2}})_{k=1,i=1}^{N,M}[\left(\mathbf{I}_N \otimes U\right) \Pi [{\bf I}_{M}\otimes (D_{\alpha}^{-1}\mathbb{F})]]^{-1}.\notag
\end{align}
And our ABAC preconditioner for \eqref{symmetrizedblttsys} is defined as:
\begin{equation}\label{precdef}
\mathcal{P}_\alpha:=\left|\mathcal{C}_\alpha\right|.
\end{equation}

 Combining Assumption \ref{akassumps} and \eqref{hatcalphaidef}, we see that $\widehat{C}_{\alpha, i}$ is a strictly diagonally dominant matrix with positive diagonal entries for $\alpha \in (0, 1]$. Therefore, $\widehat{C}_{\alpha, i}$ is clearly non-singular, and its eigenvalues lie in the right half of the complex plane. More details will be discussed on invertibility of $\widehat{C}_{\alpha, i}$ in Section \ref{convanaly}. This together with \eqref{calphadiagform} demonstrate the invertibility of our ABAC preconditioner. The invertibility is a necessary virtue for a preconditioner. Thus, throughout this paper, we set $\alpha\in(0,1]$.
 
 With the ABAC preconditioner, we employ preconditioned MINRES solver to solve the symmetric system \eqref{symmetrizedblttsys}.

\section{Convergence analysis of the proposed MINRES method}\label{convanaly}
In this section, we analyze the convergence rate of the preconditioned MINRES solver with $\mathcal{P}_{\alpha}$ as preconditioner for solving the symmetric linear system \eqref{symmetrizedblttsys}. 

\subsection{Properties of $\mathcal{C}_{\alpha}^{\frac{1}{2}}$}
In this subsection, we investigate  properties of  $\mathcal{C}_{\alpha}^{\frac{1}{2}}$ as preliminaries for the analysis in the forthcoming subsection. 

The spectrum of a square matrix $\mathbf{C}$ is denoted by $\Sigma(\mathbf{C})$. The following lemma will be used to prove that  $\mathcal{C}_{\alpha}^{\frac{1}{2}}$  is actually real-valued.

\begin{lemma}\label{lemma:3-1}
 For a diagonal matrix $D=\operatorname{diag}\left(d_i\right)_{i=1}^N \in \mathbb{C}^{N \times N}$ with $d_1 \in \mathbb{R}, d_{N-i}=\operatorname{conj}\left(d_{i+2}\right)$ for $i=0,1, \ldots, N-2$. It holds that $\mathbb{F} D \mathbb{F}^* \in \mathbb{R}^{N \times N}$.  
\end{lemma}
\begin{proof}
In order to prevent any potential confusion arising from the use of the same symbol, let us denote $\theta_N$ as $\theta$.
\begin{align*}
\left(\mathbb{F} D \mathbb{F}^*\right)(i, j)=\sum_{m=1}^N \mathbb{F}(i, m) d_m \mathbb{F}^*(m, j) & =\sum_{m=1}^N \theta^{(i-1)(m-1)} d_m \theta^{-(m-1)(j-1)} \\
& =\sum_{m=1}^N d_m \theta^{(m-1)(i-j)} \\
& =d_1+\sum_{m=0}^{N-2} d_{m+2} \theta^{(m+1)(i-j)} .
\end{align*}
If $N$ is odd, then
\begin{align*}
\left(\mathbb{F} D \mathbb{F}^* \right)(i, j) & =d_1+\sum_{m=0}^{\left\lfloor\frac{N}{2}\right\rfloor-1} \left[d_{m+2} \theta^{(m+1)(i-j)}+d_{N-m} \theta^{(N-m-1)(i-j)} \right]\\
& =d_1+\sum_{m=0}^{\left\lfloor\frac{N}{2}\right\rfloor-1} \left[ d_{m+2} \theta^{(m+1)(i-j)}+\operatorname{conj}\left(d_{m+2} \theta^{(m+1)(i-j)}\right) \right] \in \mathbb{R} .
\end{align*}
If $N$ is even, then
\begin{align*}
\left(\mathbb{F} D \mathbb{F}^*\right)(i, j)=d_1+d_{\frac{N}{2}+1} \theta^{\left(\frac{N}{2}\right)(i-j)}+\sum_{m=0}^{\frac{N}{2}-2} \left[d_{m+2} \theta^{(m+1)(i-j)}+d_{N-m} \theta^{(N-m-1)(i-j)}\right] .
\end{align*}
Notice that, when $N$ is even, $d_{\frac{N}{2}+1}=d_{N-\left(\frac{N}{2}-1\right)}=\operatorname{conj}\left(d_{\frac{N}{2}-1+2}\right)=\operatorname{conj}\left(d_{\frac{N}{2}+1}\right)$ implies that $d_{\frac{N}{2}+1} \in \mathbb{R}$ and $\theta^{\left(\frac{N}{2}\right)(i-j)}=\exp (\pi \mathbf{i}(i-j))=(-1)^{i-j} \in \mathbb{R}$. Therefore, if $N$ is even, then
$$
\begin{aligned}
\left(\mathbb{F} D \mathbb{F}^*\right)(i, j) & =d_1+(-1)^{i-j} d_{\frac{N}{2}+1}+\sum_{m=0}^{\frac{N}{2}-2} \left[ d_{m+2} \theta^{(m+1)(i-j)}+d_{N-m} \theta^{(N-m-1)(i-j)}\right]\\
& =d_1+(-1)^{i-j} d_{\frac{N}{2}+1}+\sum_{m=0}^{\frac{N}{2}-2} \left[d_{m+2} \theta^{(m+1)(i-j)}+\operatorname{conj}\left( d_{m+2} \theta^{(m+1)(i-j)}\right)\right] \in \mathbb{R}.
\end{aligned}
$$
The proof is completed.
\end{proof}

Recall that $A_{(k)}$'s ($0\leq k\leq N-1$) are simultaneously diagonalizable with
\begin{equation*}
	A_{(k)}=U{\rm diag}(\lambda_i^{(k)})U^{\rm T},\quad k=0,1,...,N-1.
\end{equation*}
By Assumption \ref{akassumps}${\bf (ii)}$, we immediately have the following proposition.
\begin{proposition}\label{assumpinterpproposi}
\begin{equation*}
\min\limits_{1\leq i\leq M}\left(\lambda_{i}^{(0)}-\sum\limits_{j=1}^{N-1}|\lambda_{i}^{(j)}|\right)\geq c_0>0.
\end{equation*}
\end{proposition}

\begin{lemma}\label{calphasqrtlm}
Let $\mathcal{C}_\alpha \in \mathbb{R}^{M N \times M N}$ be defined by \eqref{calphadiagform} with $0<\alpha \leq 1$. Then,\\
(a) $\mathcal{C}_\alpha^{\frac{1}{2}}$ is invertible;\\
(b) $\mathcal{C}_\alpha^{\frac{1}{2}}$ is a real-valued matrix;\\
(c) $\mathcal{Y} \mathcal{C}_\alpha^{\frac{1}{2}}$ is real symmetric.
\end{lemma}
\begin{proof}
By the definition of $\mathcal{C}_\alpha^{\frac{1}{2}}$, we see that the invertibility of $\mathcal{C}_\alpha^{\frac{1}{2}}$ is equivalent to invertibility of $\mathcal{C}_\alpha$. By \eqref{calphadiagform}, we see that to show the invertibility of \eqref{calphadiagform}, it suffices to show the invertibility of $\widehat{C}_{\alpha, i}$ for each $i=1,2,...,M$. By the definition of $\widehat{C}_{\alpha, i}$ given in \eqref{hatcalphaidef}, we see that
\begin{align*}
|\widehat{C}_{\alpha, i}(k,k)|-\sum\limits_{1\leq j\leq N,j\neq k} |\widehat{C}_{\alpha, i}(k,j)|\geq |\lambda_{i}^{(0)}|-\sum\limits_{j=1}^{N-1} |\lambda_{i}^{(j)}|,\quad k=1,2,...,N.
\end{align*}

By Proposition \ref{assumpinterpproposi}, we have $\lambda_i^{(0)}>0$ and for each $i=1,2,...,M$, it holds that
\begin{equation*}
|\widehat{C}_{\alpha, i}(k,k)|-\sum\limits_{1\leq j\leq N,j\neq k} |\widehat{C}_{\alpha, i}(k,j)|\geq \lambda_i^{(0)}-\sum\limits_{j=1}^{N-1}|\lambda_{i}^{(j)}|\geq c_0>0,\quad k=1,2,...,N.
\end{equation*}
That means  $\widehat{C}_{\alpha, i}$  is a strictly diagonally dominant matrix with positive diagonal entries, which proves the invertibility of $\widehat{C}_{\alpha, i}$ for each $i=1,2,...,M$. This completes the proof of part (a).

By \eqref{sqrtcalphabdiagform} and \eqref{hatcalphaidiagform}, we see that
\begin{equation}\label{calphasqrtbdiagform}
\mathcal{C}_\alpha^{\frac{1}{2}}=\left(\mathbf{I}_N \otimes U\right) \Pi \operatorname{blockdiag}\left(\widehat{C}_{\alpha,i}^{\frac{1}{2}}\right)_{i=1}^M \Pi^{\rm T}\left(\mathbf{I}_N \otimes U\right)^{\rm T},
\end{equation}
with 
\begin{equation*}
\widehat{C}_{\alpha,i}^{\frac{1}{2}}=D_\alpha^{-1} \mathbb{F} \Lambda_{\alpha, i}^{\frac{1}{2}} \mathbb{F}^* D_\alpha,\quad i=1,2,...,M.
\end{equation*}
Hence, to show $\mathcal{C}_\alpha^{\frac{1}{2}}$ is real, it suffices to show that $\widehat{C}_{\alpha, i}^{\frac{1}{2}}$ is real-valued for each $i$.
We will employ Lemma \ref{lemma:3-1} to show this. Recall that $ \Lambda_{\alpha, i}^{\frac{1}{2}}={\rm diag}\left((\lambda_i^{(k,\alpha)})^{\frac{1}{2}}\right)_{k=1}^{N}$.
By \eqref{hatcalphaieigvaldef} and the fact that $\lambda_i^{(k)}$'s as eigenvalues of Hermitian matrices are all real numbers for $0\leq k\leq N-1,1\leq i\leq M$, we have
\begin{equation*}
\begin{aligned}
    \lambda_i^{(1, \alpha)}&=\sum_{j=1}^N \widehat{C}_{\alpha, i}(j, 1) \alpha^{\frac{j-1}{N}}\\
    &=\lambda_i^{(0)}+\lambda_i^{(1)} \alpha^{\frac{1}{N}}+\lambda_i^{(2)}\alpha^{\frac{2}{N}}+\dots+\lambda_i^{(N-1)}\alpha^{\frac{N-1}{N}}\\
    &>\lambda_i^{(0)}-\sum_{j=1}^{N-1}\left|\lambda_i^{(j)}\right|\\
    &>0.
\end{aligned}
\end{equation*}
 Therefore, $0<( \lambda_i^{(1, \alpha)})^{\frac{1}{2}} \in \mathbb{R}$. Moreover,
$$
\begin{aligned}
 \lambda_i^{(N-k, \alpha)} & =\sum_{j=1}^N  \widehat{C}_{\alpha, i}(j, 1) \alpha^{\frac{j-1}{N}} \theta^{-(N-k-1)(j-1)} \\
& =\sum_{j=1}^N  \widehat{C}_{\alpha, i}(j, 1) \alpha^{\frac{j-1}{N}} \theta^{(k+1)(j-1)} \\
& =\sum_{j=1}^N  \widehat{C}_{\alpha, i}(j, 1) \alpha^{\frac{j-1}{N}} \theta^{(k+2-1)(j-1)} \\
& =\operatorname{conj}\left(\sum_{j=1}^N \widehat{C}_{\alpha, i}(j, 1) \alpha^{\frac{j-1}{N}} \theta^{-(k+2-1)(j-1)}\right) \\
& =\operatorname{conj}\left( \lambda_i^{(k+2, \alpha)}\right), \quad k=0,1, \ldots, N-2.
\end{aligned}
$$
Accordingly, we have
$$
( \lambda_i^{(N-k, \alpha)})^{\frac{1}{2}}=\operatorname{conj}\left(( \lambda_i^{(k+2, \alpha)})^{\frac{1}{2}}\right), \quad K=0,1, \ldots, N-2.
$$
By the Lemma \ref{lemma:3-1}, it is easy to obtain $\mathbb{F} \Lambda_{\alpha, i}^{\frac{1}{2}} \mathbb{F}^*  \in \mathbb{R}^{N \times N}$. Consequently, $\widehat{C}_{\alpha, i}^{\frac{1}{2}}$ is also real-valued, which implies that $\mathcal{C}_\alpha^{\frac{1}{2}}$ is real-valued. This completes the proof of part (b).

We now turn to the proof of part (c). It is straightforward to verify that $\widehat{C}_{\alpha,i}^{\frac{1}{2}}$ is an $\alpha$-circulant matrix for each $i$. Then, from \eqref{calphasqrtbdiagform}, we see that $\mathcal{C}_{\alpha}^{\frac{1}{2}}$ is actually a block $\alpha$-circulant matrix. To see this, we notice that
\begin{equation*}
 \Pi \operatorname{blockdiag}\left(\widehat{C}_{\alpha,i}^{\frac{1}{2}}\right)_{i=1}^M \Pi^{\rm T}=\left[\begin{array}{ccccc}
 	\tilde{A}_0 & \alpha 	\tilde{A}_{N-1} & \cdots & \cdots & \alpha 	\tilde{A}_{1}  \\
 	\tilde{A}_1  & \ddots & \ddots & \ddots & \vdots \\
 	\vdots & \ddots & \ddots & \ddots & \vdots \\
 	\vdots & \ddots & \ddots & \ddots & \alpha \tilde{A}_{N-1} \\
 \tilde{A}_{N-1} & \cdots & \cdots & \tilde{A}_1&\tilde{A}_0
 \end{array}\right] ,
\end{equation*}
with
\begin{equation*}
\tilde{A}_k={\rm diag}\left(\widehat{C}_{\alpha,i}^{\frac{1}{2}}(k+1,1)\right)_{i=1}^{M},\quad k=0,1,...,N-1.
\end{equation*}
Then, 
\begin{align*}
\mathcal{C}_\alpha^{\frac{1}{2}}&=\left(\mathbf{I}_N \otimes U\right) \Pi \operatorname{blockdiag}\left(\widehat{C}_{\alpha,i}^{\frac{1}{2}}\right)_{i=1}^M \Pi^{\rm T}\left(\mathbf{I}_N \otimes U\right)^{\rm T}\\
&=\left[\begin{array}{ccccc}
	\hat{A}_0 & \alpha 	\hat{A}_{N-1} & \cdots & \cdots & \alpha 	\hat{A}_{1}  \\
	\hat{A}_1  & \ddots & \ddots & \ddots & \vdots \\
	\vdots & \ddots & \ddots & \ddots & \vdots \\
	\vdots & \ddots & \ddots & \ddots & \alpha \hat{A}_{N-1} \\
	\hat{A}_{N-1} & \cdots & \cdots & \hat{A}_1&\hat{A}_0
\end{array}\right],\quad \hat{A}_k=U\tilde{A}_k U^{\rm T},\quad k=0,1,...,N-1.
\end{align*}
    Hence, $\mathcal{C}_{\alpha}^{\frac{1}{2}}$ is actually a block $\alpha$-circulant matrix with symmetric blocks. Moreover, as $\widehat{C}_{\alpha,i}^{\frac{1}{2}}$ is real-valued, $\mathcal{C}_{\alpha}^{\frac{1}{2}}$ is a block Toeplitz matrix with real symmetric blocks. That means $\mathcal{Y}\mathcal{C}_{\alpha}^{\frac{1}{2}}$ is a block Hankel matrix with real symmetric blocks, which is therefore real symmetric.

 \end{proof}

\subsection{Convergence analysis}
In this subsection, we resort to the following lemma as framework for analyzing the convergence behavior of the proposed preconditioned MINRES method.
\begin{lemma}\textnormal{\cite{elman2014finite}}\label{minrescvglm}
Let $\mathbf{P} \in \mathbb{R}^{n_0 \times n_0}$ and $\mathbf{A} \in \mathbb{R}^{n_0 \times n_0}$ be a symmetric positive definite matrix and a symmetric nonsingular matrix, respectively. Suppose $\Sigma\left(\mathbf{P}^{-1} {\bf A} \right) \in\left[-a_1,-a_2\right] \cup\left[a_3, a_4\right]$ with $a_4 \geq a_3>0, a_1 \geq a_2>0$ and $a_1-a_2=$ $a_4-a_3$. Then, the MINRES solver with $\mathbf{P}$ as a preconditioner for the linear system $\mathbf{A x}=\mathbf{y} \in \mathbb{R}^{n_0 \times 1}$ with arbitrarily given $\mathbf{y} \in \mathbb{R}^{n_0 \times 1}$ satisfies the following convergence estimation
$$
\left\|\mathbf{r}_k\right\|_2 \leqslant 2\left(\frac{\sqrt{a_1 a_4}-\sqrt{a_2 a_3}}{\sqrt{a_1 a_4}+\sqrt{a_2 a_3}}\right)^{\lfloor k / 2\rfloor}\left\|\mathbf{r}_0\right\|_2,
$$
where $\mathbf{r}_k=\mathbf{P}^{-1} \mathbf{y}-\mathbf{P}^{-1} \bf A \mathbf{x}_k$ denotes the residual vector at the $k$th iteration with $\mathbf{x}_k$ $(k \geqslant 1)$ being the $k$th iterative solution by MINRES; $\mathbf{x}_0$ denotes an arbitrary real-valued initial guess; $\lfloor k / 2\rfloor$ denotes the integer part of $k / 2$.
\end{lemma}

In order to apply Lemma \ref{minrescvglm} on analyzing convergence of the preconditioned MINRES solver with $\mathcal{P}_{\alpha}$ as preconditioner for solving the symmetric system \eqref{symmetrizedblttsys}, it is necessary to examine the spectrum of the preconditioned matrix $\mathcal{P}_\alpha^{-1} \mathcal{Y} \mathcal{A}$. Using matrix similarity, we can express this as:
$$
\Sigma\left(\mathcal{P}_\alpha^{-1} \mathcal{Y} \mathcal{A}\right)=\Sigma\left(\mathcal{P}_\alpha^{-\frac{1}{2}} \mathcal{Y} \mathcal{A} \mathcal{P}_\alpha^{-\frac{1}{2}}\right).
$$
Therefore, it is sufficient to analyze the spectral distribution of $\mathcal{P}_\alpha^{-\frac{1}{2}} \mathcal{Y} \mathcal{A} \mathcal{P}_\alpha^{-\frac{1}{2}}$. In order to do so, we decompose the preconditioned matrix $\mathcal{P}_\alpha^{-\frac{1}{2}} \mathcal{Y} \mathcal{A} \mathcal{P}_\alpha^{-\frac{1}{2}}$ as follows:
\begin{equation}\label{eq:similarpresys}
\mathcal{P}_\alpha^{-\frac{1}{2}} \mathcal{Y} \mathcal{A} \mathcal{P}_\alpha^{-\frac{1}{2}}=\mathcal{Q}_\alpha-\mathcal{E}_\alpha,
\end{equation}
where
$$
\mathcal{Q}_\alpha=\mathcal{P}_\alpha^{-\frac{1}{2}} \mathcal{Y} \mathcal{C}_\alpha \mathcal{P}_\alpha^{-\frac{1}{2}}, \quad \mathcal{E}_\alpha=\mathcal{P}_\alpha^{-\frac{1}{2}} \mathcal{Y} \mathcal{R}_\alpha \mathcal{P}_\alpha^{-\frac{1}{2}}, \quad \mathcal{R}_\alpha:=\mathcal{C}_\alpha-\mathcal{A}.
$$
The following lemma demonstrates that $\mathcal{Q}_\alpha$ is a real symmetric orthogonal matrix, meaning that $\Sigma\left(\mathcal{Q}_\alpha\right) \subset\{-1,1\}$, and that $\mathcal{E}_\alpha$ is a matrix with small norm for an appropriately selected small $\alpha$.
\begin{lemma}\label{lemma:4}
The matrix $\mathcal{Q}_\alpha=\mathcal{P}_\alpha^{-\frac{1}{2}} \mathcal{Y} \mathcal{C}_\alpha \mathcal{P}_\alpha^{-\frac{1}{2}}$ is both real symmetric and orthogonal, i.e., $\Sigma\left(\mathcal{Q}_\alpha\right) \subset\{-1,1\}$.
\end{lemma}
\begin{proof}
To commence, we verify that $\mathcal{Q}_\alpha$ is a matrix comprised of real-valued elements. Using Lemma \ref{calphasqrtlm} (b), we can deduce that both $\mathcal{C}_\alpha^{-\frac{1}{2}}$ and $\left(\mathcal{C}_\alpha^{-\frac{1}{2}}\right)^{\rm T}$ are real-valued matrices. Since $\mathcal{P}_\alpha^{-1}=\mathcal{C}_\alpha^{-\frac{1}{2}}\left(\mathcal{C}_\alpha^{-\frac{1}{2}}\right)^{\rm T}$ is also a real-valued and normal matrix, it follows that $\mathcal{P}_\alpha^{-1}$ is SPD. 
Because of the fact $\mathcal{P}_\alpha^{-1}$ is a nonsingular matrix by Lemma \ref{calphasqrtlm} (a), its matrix square root $\mathcal{P}_\alpha^{-\frac{1}{2}}$ is also a real-valued SPD matrix. As a result, $\mathcal{Q}_\alpha=\mathcal{P}_\alpha^{-\frac{1}{2}} \mathcal{Y} \mathcal{C}_\alpha \mathcal{P}_\alpha^{-\frac{1}{2}}$ is a product of real-valued matrices and is therefore a real matrix.
Using the fact that $\left(\mathcal{C}_\alpha^{\frac{1}{2}}\right)^{\rm T} \mathcal{Y}=\mathcal{Y} \mathcal{C}_\alpha^{\frac{1}{2}}$ and $\left(\mathcal{C}_\alpha^{-\frac{1}{2}}\right)^{\rm T} \mathcal{Y}=\mathcal{Y}\mathcal{C}_\alpha^{-\frac{1}{2}}$ in Lemma \ref{calphasqrtlm} (c), we find
\begin{equation*}
\begin{aligned}
\mathcal{Q}_\alpha^{\mathrm{T}} \mathcal{Q}_\alpha=\mathcal{Q}_\alpha^2 & =\mathcal{P}_\alpha^{-\frac{1}{2}} \mathcal{Y} \mathcal{C}_\alpha \mathcal{P}_\alpha^{-1} \mathcal{Y} \mathcal{C}_\alpha \mathcal{P}_\alpha^{-\frac{1}{2}} \\
& =\mathcal{P}_\alpha^{-\frac{1}{2}} \mathcal{Y} \mathcal{C}_\alpha \mathcal{C}_\alpha^{-\frac{1}{2}}\left(\mathcal{C}_\alpha^{-\frac{1}{2}}\right)^{\rm T} \mathcal{Y}\mathcal{C}_\alpha \mathcal{P}_\alpha^{-\frac{1}{2}} \\
& =\mathcal{P}_\alpha^{-\frac{1}{2}} \mathcal{Y} \mathcal{C}_\alpha^{\frac{1}{2}} \mathcal{Y} \mathcal{C}_\alpha^{\frac{1}{2}} \mathcal{P}_\alpha^{-\frac{1}{2}} \\
& =\mathcal{P}_\alpha^{-\frac{1}{2}}\left(\mathcal{C}_\alpha^{\frac{1}{2}}\right)^{\rm T} \underbrace{\mathcal{Y}^2}_{=I_{m n}} \mathcal{C}_\alpha^{\frac{1}{2}} \mathcal{P}_\alpha^{-\frac{1}{2}} \\
& =\mathcal{P}_\alpha^{-\frac{1}{2}} \mathcal{P}_\alpha \mathcal{P}_\alpha^{-\frac{1}{2}} \\
& =\mathbf{I}_{M N}.
\end{aligned}
\end{equation*}
As $\mathcal{Q}_\alpha$ is a real symmetric and orthogonal matrix, it follows that its spectrum $\Sigma\left(\mathcal{Q}_\alpha\right)$ is a subset of the set of real numbers that have absolute value equal to 1, i.e., $\Sigma\left(\mathcal{Q}_\alpha\right) \subset \mathbb{R} \cap\{z:|z|=1\}=\{-1,1\}$. Consequently, we have shown that $\mathcal{Q}_\alpha$ is a real symmetric orthogonal matrix. This completes the proof.
\end{proof}

In what follows, we will estimate $\left\|\mathcal{E}_\alpha\right\|_2$.

Based on the proof of Lemma \ref{lemma:4}, we know that $\mathcal{P}_\alpha^{-\frac{1}{2}}$ is a real symmetric matrix. Moreover, since $\mathcal{R}_\alpha$ is a real block Toeplitz matrix with symmetric blocks, it follows that $\mathcal{Y} \mathcal{R}_\alpha$ is also a real symmetric matrix. Therefore, we have $\mathcal{E}_\alpha=\mathcal{P}_\alpha^{-\frac{1}{2}} \mathcal{Y} \mathcal{R}_\alpha \mathcal{P}_\alpha^{-\frac{1}{2}}$, which is a real symmetric matrix. Here, $\rho(\cdot)$ denotes the spectral radius of a matrix. We have
\begin{equation}\label{eq:ealpha2norm}
\begin{aligned}
\left\|\mathcal{E}_\alpha\right\|_2=\left\|\mathcal{P}_\alpha^{-\frac{1}{2}} \mathcal{Y} \mathcal{R}_\alpha \mathcal{P}_\alpha^{-\frac{1}{2}}\right\|_2 & =\rho\left(\mathcal{P}_\alpha^{-\frac{1}{2}} \mathcal{Y} \mathcal{R}_\alpha \mathcal{P}_\alpha^{-\frac{1}{2}}\right) \\
& =\rho\left(\mathcal{P}_\alpha^{-1} \mathcal{Y} \mathcal{R}_\alpha\right) \\
& \leqslant\left\|\mathcal{P}_\alpha^{-1} \mathcal{Y} \mathcal{R}_\alpha\right\|_2.
\end{aligned}
\end{equation}
Due to the block orthogonal diagonalization formula \eqref{calphadiagform} for $\mathcal{C}_\alpha$, the same diagonalization also holds for $\mathcal{P}_\alpha^{-1} \mathcal{Y} \mathcal{R}_\alpha$, i.e.,
\begin{equation}\label{precdefsymmralpha}
\mathcal{P}_\alpha^{-1} \mathcal{Y} \mathcal{R}_\alpha=\left(\mathbf{I}_N \otimes U\right) \Pi \operatorname{blockdiag}\left(\widehat{C}_{\alpha, i}^{-\frac{1}{2}}\left(\widehat{C}_{\alpha, i}^{-\frac{1}{2}}\right)^{\rm T} H_{\alpha, i}\right)_{i=1}^M \Pi^{\rm T}\left(\mathbf{I}_N \otimes U\right)^{\rm T},
\end{equation}
where $\left(\mathbf{I}_N \otimes U\right) \Pi$ is an orthogonal matrix, $\widehat{C}_{\alpha, i}$ defined as \eqref{hatcalphaidef} and 
\begin{equation}\label{eq:amatdefhalpha}
H_{\alpha, i}=\left[\begin{array}{ccccc}
 &  &  &  & \\
 &  &  &  & \alpha \lambda_i^{(N-1)} \\
 &  &  & \reflectbox{$\ddots$} & \vdots \\
 &  & \reflectbox{$\ddots$} & \reflectbox{$\ddots$} & \alpha \lambda_i^{(2)} \\
 & \alpha \lambda_i^{(N-1)} & \cdots & \alpha \lambda_i^{(2)} & \alpha \lambda_i^{(1)} 

\end{array}\right] \in \mathbb{R}^{N \times N}.
\end{equation}

With the block diagonalization form \eqref{precdefsymmralpha}, we see that
\begin{equation}\label{eq:uppbpresymmsys}
    \left\|\mathcal{P}_\alpha^{-1} \mathcal{Y} \mathcal{R}_\alpha\right\|_2=\max _{1 \leqslant i \leqslant M}\left\|\widehat{C}_{\alpha, i}^{-\frac{1}{2}}\left(\widehat{C}_{\alpha, i}^{-\frac{1}{2}}\right)^{\rm T} H_{\alpha, i}\right\|_2 \leqslant \max _{1 \leqslant i \leqslant M}\left\|\widehat{C}_{\alpha, i}^{-\frac{1}{2}}\right\|_2^2\left\|H_{\alpha, i}\right\|_2.
\end{equation}
It is suggested that in order to calculate an upper bound for $\left\|\mathcal{E}_\alpha\right\|_2$, it is adequate to determine the upper bounds of $\left\|\widehat{C}_{\alpha, i}^{-\frac{1}{2}}\right\|_2$ and $\left\|H_{\alpha, i}\right\|_2$ for each $i$, individually.

In the following, we first estimate the upper bound of $\left\|H_{\alpha, i}\right\|_2$ for each $i$.
As $H_{\alpha, i}$ defined in \eqref{eq:amatdefhalpha} is clearly Hermitian matrix for each $i$, we have
\begin{equation}\label{eq:uppbhalpha}
\left\|H_{\alpha, i}\right\|_2= \rho\left(H_{\alpha, i}\right) \leqslant\left\|H_{\alpha, i}\right\|_{\infty}=\alpha\sum_{k=1}^{N-1} |\lambda_i^{(k)}|\leq \alpha(\lambda_i^{(0)}-c_0).
\end{equation}

In the next, we will estimate upper bounds of $\left\|\widehat{C}_{\alpha, i}^{-\frac{1}{2}}\right\|_2$. 

For a complex square matrix ${\bf C}$, denote
\begin{equation*}
\Re({\bf C}):=\frac{1}{2}({\bf C}+{\bf C}^{*}),\quad \Im({\bf C}):=\frac{1}{2}({\bf C}-{\bf C}^{*}).
\end{equation*}
Clearly, if ${\bf C}$ is an $1\times 1$ complex matrix (i.e., a complex number), then $\Re({\bf C})$ and $\Im({\bf C})$ are simply real and imaginary parts of ${\bf C}$, respectively.

Denote
\begin{equation*}
\mathbb{C}_{++}:=\{z\in\mathbb{C}|\Re(z)>0\},\quad \mathbb{R}_{-}:=\{x\in\mathbb{R}|x\leq 0\}.
\end{equation*}
\begin{definition}\textnormal{(see, e.g., \cite{schmitt1992})}\label{pb-squarerootdef}
For a  complex square matrix ${\bf C}$ with $\Sigma({\bf C})\subset \mathbb{C}\setminus\mathbb{R}_{-} $, there exists an unique matrix ${\bf R}$ such that ${\bf R}^2={\bf C}$ and $\Sigma({\bf R})\subset \mathbb{C}_{++}$. We call such  ${\bf R}$ as positive-branch square root of ${\bf C}$ and denote it as $\mathcal{S}({\bf C}):={\bf R}$.
\end{definition}

\begin{lemma}\textnormal{(See \cite[Lemma 2.2]{schmitt1992})}\label{sqcontilm}
Suppose ${\bf C}_1$ and ${\bf C}_2$ are two complex square matrix of same size with $\lambda_{\min}(\Re({\bf C}_i))>0$ for $i=1,2$. Then, $\mathcal{S}({\bf C}_i)$ uniquely exists for $i=1,2$ and
\begin{equation*}
\|\mathcal{S}({\bf C}_1)-\mathcal{S}({\bf C}_2)\|_2\leq \frac{1}{\sqrt{\lambda_{\min}(\Re({\bf C}_1))}+\sqrt{\lambda_{\min}(\Re({\bf C}_2))}}\|{\bf C}_1-{\bf C}_2\|_2.
\end{equation*}
\end{lemma}

Denote 
\begin{equation*}
\hat{T}_i=\left[\begin{array}[c]{cccc}
	\lambda_i^{(0)} &  &  &     \\
	\lambda_i^{(1)} & \ddots &    &   \\
	\vdots & \ddots & \ddots &   \\
	\lambda_i^{(N-1)} & \cdots  & \lambda_i^{(1)} & \lambda_i^{(0)}
\end{array}\right]  \in \mathbb{R}^{N \times N},\quad i=1,2,...,M.
\end{equation*}

\begin{lemma}\label{sqrtcalphaapplm}
	Take $\alpha\in(0,1]$. For each $i=1,2,...,M$, $\mathcal{S}(\hat{T}_i)$ exists and
	\begin{equation*}
	\max\limits_{1\leq i\leq M}\|\mathcal{S}(\hat{T}_i)-\widehat{C}_{\alpha,i}^{\frac{1}{2}}\|_2\leq \frac{\alpha(\|A_{(0)}\|_2-c_0)}{2\sqrt{c_0}},
	\end{equation*} 
	where $c_0$ is given in Proposition \ref{assumpinterpproposi}.
\end{lemma}
\begin{proof}
By Proposition \ref{assumpinterpproposi}, it is easy to see that $\Re(\hat{T}_i)$ and $\Re(\widehat{C}_{\alpha,i})$ are both strictly diagonally dominant matrices with positive diagonal entries. Moreover, applying the well-known Gershgorin circle Theorem, it is straightforward to see that
\begin{align*}
&\lambda_{\min}(\Re(\hat{T}_i))\geq \lambda_i^{(0)}-\sum\limits_{j=0}^{N-1}|\lambda_i^{(j)}|\geq c_0>0,\\
&\lambda_{\min}(\Re(\widehat{C}_{\alpha,i}))\geq \lambda_i^{(0)}-\sum\limits_{j=0}^{N-1}|\lambda_i^{(j)}|\geq c_0>0.
\end{align*}
Then, from Lemma \ref{sqcontilm}, we see that both $\mathcal{S}(\hat{T}_i)$ and $\mathcal{S}(\widehat{C}_{\alpha,i})$ exist and
\begin{align*}
\|\mathcal{S}(\hat{T}_i)-\mathcal{S}(\widehat{C}_{\alpha,i})\|_2&\leq \frac{1}{2\sqrt{c_0}}\|\hat{T}_i-\widehat{C}_{\alpha,i}\|_2\\
&\leq \frac{1}{2\sqrt{c_0}}\sqrt{\|\hat{T}_i-\widehat{C}_{\alpha,i}\|_{\infty}\|(\hat{T}_i-\widehat{C}_{\alpha,i})^{\rm T}\|_{\infty}}\\
&=\frac{\alpha}{2\sqrt{c_0}}\sum\limits_{j=0}^{N-1}|\lambda_i^{(j)}|\leq \frac{\alpha(\lambda_i^{(0)}-c_0)}{2\sqrt{c_0}}.
\end{align*}
 Since $\Sigma(\widehat{C}_{\alpha,i})\subset\mathbb{C}_{++}$ and $\widehat{C}_{\alpha,i}$ is diagonalizable, $\widehat{C}_{\alpha,i}^{\frac{1}{2}}$ exists and $\Sigma(\widehat{C}_{\alpha,i}^{\frac{1}{2}})\subset\mathbb{C}_{++}$ for each $k$ according to Definition \ref{diagonalizablematsqrtdef}.  Then, according to Definition \ref{pb-squarerootdef}, we see that $\mathcal{S}(\widehat{C}_{\alpha,i})=\widehat{C}_{\alpha,i}^{\frac{1}{2}}$. That means,
\begin{equation*}
\|\mathcal{S}(\hat{T}_i)-\widehat{C}_{\alpha,i}^{\frac{1}{2}}\|_2\leq \frac{\alpha(\lambda_i^{(0)}-c_0)}{2\sqrt{c_0}}\leq \frac{\alpha(\|A_{(0)}\|_2-c_0)}{2\sqrt{c_0}},
\end{equation*}
for each $i=1,2,...,M$. The proof is complete.
\end{proof}

\begin{lemma}\label{lemma:10}\textnormal{\cite{horn2012matrix}} For $N \times N$ Hermitian matrices $\mathbf{H}$ and $\mathbf{E}$, if $\hat{\sigma}$ is an eigenvalue of $\mathbf{H}+\mathbf{E}$, then there exists an eigenvalue $\sigma$ of $\mathbf{H}$ such that
	$$
	|\sigma-\hat{\sigma}| \leqslant\|\mathbf{E}\|_2 .
	$$
\end{lemma}
The following lemma indicates that the spectrum of the preconditioned matrix $\mathcal{P}_\alpha^{-1} \mathcal{Y} \mathcal{A}$ is located in a disjoint interval excluding the origin.

\begin{lemma}\label{lemma:11}
Take $\alpha\in(0,\nu]$ with
\begin{equation*}
\nu=\min\left\{1,\frac{2c_0\lambda_{\min}(\mathcal{S}(\mathcal{A})^{\rm T}\mathcal{S}(\mathcal{A}))}{4\sqrt{c_0}(\|A_{(0)}\|_2-c_0)\|\mathcal{S}(\mathcal{A})\|_2+(\|A_{(0)}\|_2-c_0)^2}\right\}>0.
\end{equation*}
Then,
$$
\Sigma\left(\mathcal{P}_\alpha^{-1} \mathcal{Y} \mathcal{A}\right) \subset\left[-1-\alpha\mu ,-1+ \alpha \mu\right] \cup\left[1-\alpha\mu, 1+ \alpha \mu\right],
$$
where
\begin{equation*}
\mu=\max\limits_{1\leq i\leq M}  \frac{2(\|A_{(0)}\|_2-c_0)}{\lambda_{\min}(\mathcal{S}(\mathcal{A})^{\rm T}\mathcal{S}(\mathcal{A}))}
\end{equation*}
is independent of $\alpha$.
\end{lemma}
\begin{proof}
Since $\mathcal{A}$ is nonsingular, $\mathcal{S}(\mathcal{A})$ is clearly nonsingular. Hence,\\ $\lambda_{\min}(\mathcal{S}(\mathcal{A})^{\rm T}\mathcal{S}(\mathcal{A}))>0$. That means the positive numbers $\nu$ and $\mu$ are well-defined. 

By matrix similarity and \eqref{eq:similarpresys}, we have $\Sigma\left(\mathcal{P}_\alpha^{-1} \mathcal{Y} \mathcal{A}\right)=\Sigma\left(\mathcal{P}_\alpha^{-\frac{1}{2}} \mathcal{Y} \mathcal{A} \mathcal{P}_\alpha^{-\frac{1}{2}}\right)=$ $\Sigma\left(\mathcal{Q}_\alpha-\mathcal{E}_\alpha\right)$. Since $\mathcal{Q}_\alpha$ and $-\mathcal{E}_\alpha$ are both Hermitian, Lemma \ref{lemma:10} is applicable. In other words, for $\hat{\sigma} \in \Sigma\left(\mathcal{Q}_\alpha-\mathcal{E}_\alpha\right)$, there exists $\sigma \in \Sigma\left(\mathcal{Q}_\alpha\right)$ such that $|\hat{\sigma}-\sigma| \leqslant\left\|\mathcal{E}_\alpha\right\|_2$. By Lemma \ref{lemma:4}, we know that either $\sigma=1$ or $\sigma=-1$. That means either $|\hat{\sigma}-1| \leqslant\left\|\mathcal{E}_\alpha\right\|_2$ or $|\hat{\sigma}+1| \leqslant\left\|\mathcal{E}_\alpha\right\|_2$. Thus, we have
\begin{equation*}
	\begin{aligned}
		\Sigma\left(\mathcal{P}_\alpha^{-1} \mathcal{Y} \mathcal{A}\right) & =\Sigma\left(\mathcal{Q}_\alpha-\mathcal{E}_\alpha\right) \\
		& \subset\left[-1-\left\|\mathcal{E}_\alpha\right\|_2,-1+\left\|\mathcal{E}_\alpha\right\|_2\right] \cup\left[1-\left\|\mathcal{E}_\alpha\right\|_2, 1+\left\|\mathcal{E}_\alpha\right\|_2\right].
	\end{aligned}
\end{equation*}

From \eqref{eq:ealpha2norm}, \eqref{eq:uppbpresymmsys} and \eqref{eq:uppbhalpha}, we see that
\begin{align}
\|\mathcal{E}_{\alpha}\|_2\leq \max\limits_{1\leq i\leq M} \left\|\widehat{C}_{\alpha, i}^{-\frac{1}{2}}\right\|_2^2\left\|H_{\alpha, i}\right\|_2
&\leq \alpha \max\limits_{1\leq i\leq M} \left\|\widehat{C}_{\alpha, i}^{-\frac{1}{2}}\right\|_2^2(\lambda_i^{(0)}-c_0)\notag\\
&=\alpha \max\limits_{1\leq i\leq M} \frac{\lambda_i^{(0)}-c_0}{\lambda_{\min}\left(\widehat{C}_{\alpha, i}^{\frac{{\rm T}}{2}}\widehat{C}_{\alpha, i}^{\frac{1}{2}}\right)}.\label{ealphaesti1}
\end{align}
For each $i$, it is straightforward to verify that
\begin{align*}
\widehat{C}_{\alpha, i}^{\frac{{\rm T}}{2}}\widehat{C}_{\alpha, i}^{\frac{1}{2}}=&\mathcal{S}(\hat{T}_i)^{\rm T}\mathcal{S}(\hat{T}_i)+\mathcal{S}(\hat{T}_i)^{\rm T}(\widehat{C}_{\alpha, i}^{\frac{1}{2}}-\mathcal{S}(\hat{T}_i))+(\widehat{C}_{\alpha, i}^{\frac{1}{2}}-\mathcal{S}(\hat{T}_i))^{\rm T}\mathcal{S}(\hat{T}_i)\\
&+(\widehat{C}_{\alpha, i}^{\frac{1}{2}}-\mathcal{S}(\hat{T}_i))^{\rm T}(\widehat{C}_{\alpha, i}^{\frac{1}{2}}-\mathcal{S}(\hat{T}_i)).
\end{align*}
Therefore, for each $i$, we have
\begin{align*}
\lambda_{\min}(\widehat{C}_{\alpha, i}^{\frac{{\rm T}}{2}}\widehat{C}_{\alpha, i}^{\frac{1}{2}})\geq& \lambda_{\min}(\mathcal{S}(\hat{T}_i)^{\rm T}\mathcal{S}(\hat{T}_i))-\|\mathcal{S}(\hat{T}_i)^{\rm T}\|_2\|\widehat{C}_{\alpha, i}^{\frac{1}{2}}-\mathcal{S}(\hat{T}_i)\|_2\\
&-\|(\widehat{C}_{\alpha, i}^{\frac{1}{2}}-\mathcal{S}(\hat{T}_i))^{\rm T}\|_2\|\mathcal{S}(\hat{T}_i)\|_2-\|\widehat{C}_{\alpha, i}^{\frac{1}{2}}-\mathcal{S}(\hat{T}_i)\|_2^2\\
=&\lambda_{\min}(\mathcal{S}(\hat{T}_i)^{\rm T}\mathcal{S}(\hat{T}_i))-2\|\mathcal{S}(\hat{T}_i)\|_2\|\widehat{C}_{\alpha, i}^{\frac{1}{2}}-\mathcal{S}(\hat{T}_i)\|_2-\|\widehat{C}_{\alpha, i}^{\frac{1}{2}}-\mathcal{S}(\hat{T}_i)\|_2^2.
\end{align*}
Since $\alpha\in(0,\nu]\subset(0,1]$, Lemma \ref{sqrtcalphaapplm} implies that
\begin{align*}
\lambda_{\min}(\widehat{C}_{\alpha, i}^{\frac{{\rm T}}{2}}\widehat{C}_{\alpha, i}^{\frac{1}{2}})& \geq \lambda_{\min}(\mathcal{S}(\hat{T}_i)^{\rm T}\mathcal{S}(\hat{T}_i))-\frac{\alpha(\|A_{(0)}\|_2-c_0)\|\mathcal{S}(\hat{T}_i)\|_2}{\sqrt{c_0}}-\frac{\alpha^2(\|A_{(0)}\|_2-c_0)^2}{4c_0}\\
&\geq \lambda_{\min}(\mathcal{S}(\hat{T}_i)^{\rm T}\mathcal{S}(\hat{T}_i))-\frac{\alpha(\|A_{(0)}\|_2-c_0)\|\mathcal{S}(\hat{T}_i)\|_2}{\sqrt{c_0}}-\frac{\alpha(\|A_{(0)}\|_2-c_0)^2}{4c_0},
\end{align*}
for each $i$. Note that
\begin{align*}
 \|\mathcal{S}(\mathcal{A})\|_2=\max\limits_{1\leq i\leq M}\|\mathcal{S}(\hat{T}_i)\|_2,\quad \lambda_{\min}(\mathcal{S}(\mathcal{A})^{\rm T}\mathcal{S}(\mathcal{A}))=\min\limits_{1\leq i\leq M}\lambda_{\min}(\mathcal{S}(\hat{T}_i)^{\rm T}\mathcal{S}(\hat{T}_i)).
\end{align*}
Then, we have
\begin{equation*}
\min\limits_{1\leq i\leq M}\lambda_{\min}(\widehat{C}_{\alpha, i}^{\frac{{\rm T}}{2}}\widehat{C}_{\alpha, i}^{\frac{1}{2}})\geq  \lambda_{\min}(\mathcal{S}(\mathcal{A})^{\rm T}\mathcal{S}(\mathcal{A}))-\frac{\alpha(\|A_{(0)}\|_2-c_0)\|\mathcal{S}(\mathcal{A})\|_2}{\sqrt{c_0}}-\frac{\alpha(\|A_{(0)}\|_2-c_0)^2}{4c_0}.
\end{equation*}
Since $$0<\alpha\leq \nu\leq \frac{2c_0\lambda_{\min}(\mathcal{S}(\mathcal{A})^{\rm T}\mathcal{S}(\mathcal{A}))}{4\sqrt{c_0}(\|A_{(0)}\|_2-c_0)\|\mathcal{S}(\mathcal{A})\|_2+(\|A_{(0)}\|_2-c_0)^2},$$ it is easy to verify that
\begin{equation*}
\lambda_{\min}(\mathcal{S}(\mathcal{A})^{\rm T}\mathcal{S}(\mathcal{A}))-\frac{\alpha(\|A_{(0)}\|_2-c_0)\|\mathcal{S}(\mathcal{A})\|_2}{\sqrt{c_0}}-\frac{\alpha(\|A_{(0)}\|_2-c_0)^2}{4c_0}\geq \frac{1}{2}\lambda_{\min}(\mathcal{S}(\mathcal{A})^{\rm T}\mathcal{S}(\mathcal{A})),
\end{equation*}
for each $i$. Therefore,
\begin{equation*}
\lambda_{\min}(\widehat{C}_{\alpha, i}^{\frac{{\rm T}}{2}}\widehat{C}_{\alpha, i}^{\frac{1}{2}})\geq  \frac{1}{2}\lambda_{\min}(\mathcal{S}(\mathcal{A})^{\rm T}\mathcal{S}(\mathcal{A}))
\end{equation*}
for each $i$. Combining these inequalities with \eqref{ealphaesti1}, we have
\begin{equation*}
\|\mathcal{E}_{\alpha}\|_2\leq \alpha\max\limits_{1\leq i\leq M}\frac{2(\lambda_i^{(0)}-c_0)}{\lambda_{\min}(\mathcal{S}(\mathcal{A})^{\rm T}\mathcal{S}(\mathcal{A}))}\leq \alpha\max\limits_{1\leq i\leq M}\frac{2(\|A_{(0)}\|_2-c_0)}{\lambda_{\min}(\mathcal{S}(\mathcal{A})^{\rm T}\mathcal{S}(\mathcal{A}))}=\alpha\mu.
\end{equation*}
Hence,
\begin{align*}
	\Sigma\left(\mathcal{P}_\alpha^{-1} \mathcal{Y} \mathcal{A}\right)&\subset\left[-1-\left\|\mathcal{E}_\alpha\right\|_2,-1+\left\|\mathcal{E}_\alpha\right\|_2\right] \cup\left[1-\left\|\mathcal{E}_\alpha\right\|_2, 1+\left\|\mathcal{E}_\alpha\right\|_2\right]\\
	&\subset \left[-1-\alpha\mu ,-1+ \alpha \mu\right] \cup\left[1-\alpha\mu, 1+ \alpha \mu\right].
\end{align*}
The proof is complete.
\end{proof}

With Lemmas \ref{minrescvglm} and \ref{lemma:11}, one can immediately obtain the following theorem.
\begin{theorem}\label{mainthm}
For any constant $\delta \in(0,1)$, choose $\alpha \in\left(0, \zeta \right]$ with
$$
\zeta:=\min \left\{\frac{\delta^2}{\mu}, \nu\right\} ,
$$
with $\mu$ and $\nu$ defined in Lemma \ref{lemma:11}.
Then, the MINRES solver for the preconditioned system \eqref{symmetrizedblttsys} has a convergence rate independent of $M$ and $N$, i.e.,
$$
\left\|\hat{\mathbf{r}}_k\right\|_2 \leqslant 2 \delta^{k-1}\left\|\hat{\mathbf{r}}_0\right\|_2, \quad k \geqslant 1 ,
$$
where $\hat{\mathbf{r}}_k=\mathcal{P}_\alpha^{-1} \mathcal{Y} \mathcal{A} \mathbf{u}_k-\mathcal{P}_\alpha^{-1} \mathcal{Y} \mathbf{f}$ denotes the residual vector at $k$th MINRES iteration with $\mathbf{u}_k$ denoting the kth iterative solution by MINRES; $\hat{\mathbf{r}}_0$ denotes the initial residual vector computed by an arbitrary real-valued initial guess $\mathbf{u}_0$.
\end{theorem}
\begin{proof}
    Applying Lemma \ref{minrescvglm} to Lemma \ref{lemma:11}, we have $$
    \begin{aligned}
\left\|\hat{\mathbf{r}}_k\right\|_2 &\leqslant 2[\alpha\mu]^{\lfloor k / 2\rfloor}\left\|\mathbf{r}_0\right\|_2 \\
    &\leqslant 2\left[\sqrt{\alpha\mu}\right]^{k-1}\left\|\mathbf{r}_0\right\|_2.
    \end{aligned}
    $$
    Since  $\alpha \in\left(0, \zeta\right]$, we then have
    \begin{equation*}
    \sqrt{\alpha\mu}\leq \delta.
    \end{equation*}
    And hence,
    \begin{equation*}
    \left\|\hat{\mathbf{r}}_k\right\|_2 \leq  2\left[\sqrt{\alpha\mu}\right]^{k-1}\left\|\mathbf{r}_0\right\|_2\leq 2 \delta^{k-1}\left\|\mathbf{r}_0\right\|_2.
    \end{equation*}
The proof is complete.
\end{proof}
\begin{remark}
Theorem \ref{mainthm} shows the convergence rate of the proposed preconditioned MINRES solver does not deteriorate as the grid gets refined provided that $\alpha$ is sufficiently small. The theoretical results demonstrate the robustness of the proposed preconditioning technique.
\end{remark}
\section{Numerical Experiments}\label{numexams}
In this section, we firstly introduce implementation details of matrix-vector multiplication associated with the preconditioned matrix  $\mathcal{P}_\alpha^{-1} \mathcal{Y} \mathcal{A}$, which is necessary for code reproduction of the proposed solver. And then, numerical results of proposed solver on several examples are presented.

\subsection{Implementations}
When applying $\mathcal{P}_\alpha$ as preconditioner of preconditioned MINRES solver for solving \eqref{symmetrizedblttsys}, the dominant operation cost in each iteration is devoted to computing  matrix-vector product $\mathcal{P}_\alpha^{-1} \mathcal{Y} \mathcal{A}{\bf v}$ for some given vector ${\bf v}$. For fast implementation, one can compute $\mathcal{Y} \mathcal{A}{\bf v}$ first and compute $\mathcal{P}_\alpha^{-1}{\bf y}$ with ${\bf y}=\mathcal{Y} \mathcal{A}{\bf v}$ then. The blocks  $A_{(k)}$ ($0\leq k\leq N-1$) appearing in the examples tested in this section are all sparse mass or stiffness matrices from discretization of spatial operators. Exploiting sparsity of $A_{(k)}$'s and BLTT structure of  $\mathcal{A}$, the matrix-vector multiplication $\mathcal{A}{\bf v}$ can be fast computed within $\mathcal{O}(MN\log(N))$ operations by means of fast Fourier transforms (FFTs); see, e.g., \cite{ng2004iterative}. $\mathcal{Y} \mathcal{A}{\bf v}$ can be regarded as a trivial temporal reversing operation on $\mathcal{A}{\bf v}$, which cost only $\mathcal{O}(NM)$ read-write operations. 

It thus remains to discuss the computation of  $\mathcal{P}_\alpha^{-1}{\bf y}$ for a given vector ${\bf y}$. Clearly, $\mathcal{P}_\alpha^{-1} \mathbf{y}=\mathcal{C}_\alpha^{-\frac{1}{2}}\left(\mathcal{C}_\alpha^{-\frac{1}{2}}\right)^* \mathbf{y}$ for a given vector $\mathbf{y}$. We firstly discuss the implementation of computing $\left(\mathcal{C}_\alpha^{-\frac{1}{2}}\right)^* \mathbf{y}$. By \eqref{sqrtcalphabdiagform}, we have
\begin{equation*}
	\left(\mathcal{C}_\alpha^{-\frac{1}{2}}\right)^*=\left(\mathbf{I}_N \otimes U\right) \Pi \operatorname{blockdiag}\left(D_\alpha \mathbb{F} \bar{\Lambda}_{\alpha, i}^{-\frac{1}{2}} \mathbb{F}^* D_\alpha^{-1}\right)_{i=1}^M \Pi^{\rm T}\left(\mathbf{I}_N \otimes U\right)^{\rm T}.
\end{equation*}
Hence, the product $\left(\mathcal{C}_\alpha^{-\frac{1}{2}}\right)^* \mathbf{y}$ for a given vector $\mathbf{y}$ can be implemented via the following three steps:\\
	\textbf{ Step 1:} Compute $\tilde{\mathbf{y}}=\Pi^{\rm T}\left(\mathbf{I}_N \otimes U\right)^{\rm T} \mathbf{y}$;\\
	\textbf{ Step 2:} Compute $\tilde{\mathbf{z}}=\operatorname{blockdiag}\left(D_\alpha \mathbb{F} \bar{\Lambda}_{\alpha, i}^{-\frac{1}{2}} \mathbb{F}^* D_\alpha^{-1}\right)_{i=1}^M \tilde{\mathbf{y}}$;\\
	\textbf{Step 3:} Compute $\mathbf{z}=\left(\mathbf{I}_N \otimes U\right) \Pi \tilde{\mathbf{z}}$.\\
Notice that $\Pi$ and   $\Pi^{\rm T}$ are permutation matrices whose associated matrix-vector products only involve $\mathcal{O}(MN)$ read-write operations, which is fast. Moreover, for examples tested this section, ${\bf A}_k$'s are all diagonalizable by the discrete sine transform matrix $U$, which means \textbf{Step 1} and \textbf{Step 3} can be fast computed within $\mathcal{O}(MN\log M)$ operations by means of fast sine transforms (FSTs). 

Once  $\hat{\bf y}=\left(\mathcal{C}_\alpha^{-\frac{1}{2}}\right)^* \mathbf{y}$ is computed, $\mathcal{P}_\alpha^{-1} \mathbf{y}=\mathcal{C}_\alpha^{-\frac{1}{2}}\hat{\bf y}$ can be computed in a similar fashion within $\mathcal{O}(MN\log MN)$ operations.

\subsection{Numerical Results}
We test the proposed preconditioner via several examples and present the numerical results in this subsection. The experiments were conducted using GNU Octave 8.2.0 on a Dell R640 server with dual Xeon Gold 6246R 16-Cores 3.4 GHz CPUs and 512GB RAM running Ubuntu $20.04$ LTS. The CPU time in seconds was measured using the built-in functions tic/toc in Octave. All Krylov subspace solvers used in the experiments were implemented using the built-in functions in Octave. A zero initial guess is chosen for all Krylov subspace solvers unless otherwise indicated.

As previously stated in the Introduction, the BLTT system emerges from the time-space discretization of many evolutionary equations. To demonstrate the performance of our proposed ABAC preconditioner, we apply the preconditioned MINRES with ABAC preconditioner to solve the BLTT system arising from heat equations and non-local evolutionary equations, and compare its performance with the other two solvers: the preconditioned MINRES solver with the absolute-value block circulant preconditioner \cite{mcdonald2018preconditioning} and unpreconditioned MINRES solver. For ease of statement, we denote by MINRES-$\mathcal{P}_\alpha$,  MINRES-$\mathcal{P}_1$ and MINRES-$I$, the preconditioned MINRES with ABAC preconditioner, the preconditioned MINRES with absolute-value block circulant preconditioner and unpreconditioned MINRES solver, respectively.

For all the numerical experiments, we set $\alpha$  to be $10^{-8}.$ The stopping criterion for the MINRES solvers is set to be $\left\|\mathbf{r}_k\right\|_2 \leqslant 10^{-6}$, where $\mathbf{r}_k$ represents the residual vector at the $k$th iteration step. Let us denote the iteration numbers for solving a linear system using an iterative solver by `Iter'. When `Iter' $>$ 1000, we use the symbol `-'. The computational time in seconds can be denoted as `CPU', while the degree of freedom, i.e., the number of unknowns, can be denoted as `DoF'. 

For all numerical examples presented in this subsection, we consider the 2D physical square domain of form $\Omega=\left(\check{c}, \hat{c}\right) \times \left(\check{c}, \hat{c}\right)$ and it is discretized by uniform square grid with $h=\frac{\hat{c}-\check{c}}{m+1}$ being the spatial step-size along each spatial direction.  Denote $M=m^2$. For the Laplacian operator $\nabla(a({\bf x})\nabla)$ ($a({\bf x})>0$) defined on $\Omega$ appearing in the examples of this section, we adopt the central difference scheme for discretization on the uniform square grid and the resulting matrix is denoted as $\Delta_{a,h}\in\mathbb{R}^{M\times M}$. Especially, when $a\equiv 1$ (i.e., $\nabla(a({\bf x})\nabla)=\Delta$), the notation  $\Delta_{a,h}$ is simplified as $\Delta_{h}$. Clearly,  $\Delta_{a,h}$ is a real symmetric negative definite (RSND) matrix.
 Besides, we define the temporal step-size as $\tau=T/N$. Let $t_n=n \tau (0 \leq n \leq N)$ be the temporal grid points, where $N$ is the total number of time steps.

\subsection{Heat equations}
The following two examples involve the discretization of both constant coefficient and variable coefficient heat equations using the backward difference scheme (BDF) and the Crank-Nicolson scheme (CN).
\begin{equation}\label{eqn:heat var coeff}
\left\{
\begin{array}{lc}
u_{t}({\bf x},t) = \nabla(a({\bf x})\nabla u({\bf x},t)) + f({\bf x},t), & (x,t)\in \Omega \times (0,T], \\
u = 0, & ({\bf x},t)\in \partial \Omega \times (0,T], \\
u({\bf x},0)=u_0({\bf x}), & x \in \Omega\subset\mathbb{R}^{d},
\end{array}
\right.\
\end{equation}
where $\Omega$ is open, $\partial \Omega$ denotes boundary of $\Omega$;  $a, f, u_0$ are all given functions.			
In what follows, the solvers are tested for solving heat equations.

\begin{example}\label{exampleheateq}
{\rm
In this example, we consider the heat equation \eqref{eqn:heat var coeff} with constant coefficients.
\begin{align*}
&d=2,~\Omega=(0,1)\times (0,1),~T=1,~u_0({ x})=x_1(x_1-1)x_2(x_2-1),~a({\bf x})\equiv 1,\\
&f({ x},t)=\exp(t)[x_1(x_1-1)x_2(x_2-1)-2\times 10^{-6}\times(x_1(x_1-1)+x_2(x_2-1))].
\end{align*}
The exact solution of Example \ref{exampleheateq} is given by 
\begin{equation*}
u({ x},t)=\exp(t)x_1(x_1-1)x_2(x_2-1).
\end{equation*}

We adopt two discretization schemes, namely, BDF discretization and CN discretization. 
The time-space linear systems arising from the two schemes are presented in \eqref{eq:bdfheat_tssystem} \& \eqref{eq:cnheat_tssystem}, respectively.
				
The BDF discretization of \eqref{eqn:heat var coeff} gives
\begin{equation}\label{eq:bdfheat_tssystem}
\frac{1}{\tau}\begin{bmatrix}
\widetilde{L} &   &  & \\
            -I_{M}  & 	\widetilde{L}    & & \\
            & \ddots & \ddots &  \\
            &    & -I_{M} & 	\widetilde{L}
        \end{bmatrix}
        \left[\begin{array}{c}\mathbf{u}^{(1)}\\ \mathbf{u}^{(2)}\\ \vdots \\ \mathbf{u}^{(N)}
        \end{array} \right]
        =  \left[\begin{array}{c}
            \mathbf{f}^{(1)} -  {\bf u}_0 /\tau\\
            \mathbf{f}^{(2)} \\
            \vdots \\
            \mathbf{f}^{(N)}
        \end{array} \right],
    \end{equation}
    where $\widetilde{L}=I_M-\tau\Delta_h$; ${\bf f}^{(n)}$ denotes the values of $f(\cdot,t_n)$ on the spatial grid points; ${\bf u}_0$ denotes the values of $u_0$ on the spatial grid points; the unknowns of the linear system $\mathbf{u}^{(n)}$ is an approximation of $u(\cdot,t_n)$ on the spatial grid for each $n=1,2,...,N$. 
    
    The CN discretization of \eqref{eqn:heat var coeff} gives
    \begin{equation}\label{eq:cnheat_tssystem}
      \frac{1}{ \tau }  \begin{bmatrix}
            
            \widecheck{L}_1 &   &  & \\
            \widecheck{L}_2  & 	\widecheck{L}_1  & & \\
            & \ddots & \ddots &  \\
            &    & \widecheck{L}_2 & 	\widecheck{L}_1
        \end{bmatrix}
        \left[\begin{array}{c}\mathbf{u}^{(1)}\\ \mathbf{u}^{(2)}\\ \vdots \\ \mathbf{u}^{(N)}
        \end{array} \right]
        =\left[\begin{array}{c}
            \mathbf{f}^{(\frac{1}{2})} -  \frac{1}{\tau}\widecheck{L}_2{\bf u}_0 \\
            \mathbf{f}^{(\frac{3}{2})} \\
            \vdots \\
            \mathbf{f}^{(N-\frac{1}{2})}
        \end{array} \right],
    \end{equation}
    where $\widecheck{L}_1=I-\frac{\tau}{2}\Delta_h$; $\widecheck{L}_2=-I-\frac{\tau}{2}\Delta_h$; ${\bf f}^{(n-\frac{1}{2})}$ denotes the values of $f(\cdot,(n-1/2)\tau)$ on the spatial grid points;the unknowns of the linear system $\mathbf{u}^{(n)}$ is an approximation of $u(\cdot,t_n)$ on the spatial grid for each $n=1,2,...,N$. 
    
    The real symmetric blocks of the BLTT matrices given in \eqref{eq:bdfheat_tssystem} and  \eqref{eq:cnheat_tssystem}  clearly satisfy Assumption \ref{akassumps},  due to the fact that $\Delta_h$ is RSND. The corresponding proof is omitted.
We employ MINRES-$\mathcal{P}_{\alpha}$, MINRES-$\mathcal{P}_{1}$ and MINRES-$I$ to solve the all-at-once linear system arising from the BDF and CN time-space discretization of the heat equation, respectively. The corresponding numerical results are presented in Tables \ref{explbdfhteqtbl_heat_constant} \& \ref{explcnhteqtbl_constant_heat}. The two tables show that (i) MINRES-$\mathcal{P}_{\alpha}$ is more efficient than MINRES-$\mathcal{P}_{1}$ in terms of computational time and iteration number; (ii) the iteration number of MINRES-$\mathcal{P}_{\alpha}$ remains stable at 2 as the temporal and spatial grid become more refined. The stable iteration number of  MINRES-$\mathcal{P}_{\alpha}$ means that its convergence rate does not deteriorate as the grids getting refined, which supports the theoretical results presented in Theorem \ref{mainthm} and demonstrates the robustness of the proposed solver. 

\begin{table}[H]
\begin{center}
\caption{Iteration numbers and CPU times of MINRES-$\mathcal{P}_{\alpha}$, MINRES-$\mathcal{P}_1$, and MINRES-$I$ with BDF discretization for Example \ref{exampleheateq}.}\label{explbdfhteqtbl_heat_constant}
\setlength{\tabcolsep}{0.8em}
{\small
    \begin{tabular}{ccc|cc|cc|cc}
        \hline
        \multirow{2}{*}{$N$}   & \multirow{2}{*}{$m+1$} & \multirow{2}{*}{DoF} & \multicolumn{2}{c|}{MINRES-$\mathcal{P}_{\alpha}$} & \multicolumn{2}{c|}{MINRES-$\mathcal{P}_{1}$} & \multicolumn{2}{c}{MINRES-$I$} \\ \cline{4-9} 
        &                      &                      & Iter                    & CPU(s)                   & Iter                 & CPU(s)                 & Iter            & CPU(s)            \\ \hline
        \multirow{4}{*}{$2^{5}$} & $2^{5}$             & 30752                & 2                       & 0.023                    & 65                   & 0.41                   & 48              & 0.10              \\
        & $2^{6}$             & 127008               & 2                       & 0.040                    & 95                   & 1.12                   & 49              & 0.31              \\
        & $2^{7}$             & 516128               & 2                       & 0.14                     & 141                  & 5.57                   & 49              & 0.76              \\
        & $2^{8}$             & 2080800              & 2                       & 0.84                     & 185                  & 44.87                  & 50              & 3.53              \\ \hline
        \multirow{4}{*}{$2^{6}$} & $2^{5}$             & 61504                & 2                       &  0.027                    & 65                   & 0.47                   & 72              & 0.32              \\
        & $2^{6}$             & 254016               & 2                       & 0.078                    & 94                   & 2.21                   & 74              & 0.62              \\
        & $2^{7}$             & 1032256              & 2                       & 0.33                     & 140                  & 12.49                  & 75              & 2.25              \\
        & $2^{8}$             & 4161600              & 2                       & 1.99                     & 187                  & 104.01                 & 77              & 11.89             \\ \hline
        \multirow{4}{*}{$2^{7}$} & $2^{5}$             & 123008               & 2                       & 0.051                    & 65                   & 0.78                   & 133             & 0.73              \\
        & $2^{6}$             & 508032               & 2                       & 0.14                     & 94                   & 4.01                   & 133             & 1.85              \\
        & $2^{7}$             & 2064512              & 2                       & 0.73                     & 139                  & 34.35                  & 134             & 8.63              \\
        & $2^{8}$             & 8323200              & 2                       & 4.29                     & 227                  & 302.47                 & 134             & 64.70             \\ \hline
        \multirow{4}{*}{$2^{8}$} & $2^{5}$             & 246016               & 2                       & 0.076                    & 64                   & 1.26                   & 261             & 2.12              \\
        & $2^{6}$             & 1016064              & 2                       &  0.30                     & 93                   & 7.27                   & 261             & 8.50              \\
        & $2^{7}$             & 4129024              & 2                       & 1.91                     & 140                  & 73.65                  & 261             & 39.05             \\
        & $2^{8}$             & 16646400             & 2                       & 8.52                     & 227                  & 577.53                 & 261             & 250.41            \\ \hline
    \end{tabular}
}
\end{center}
\end{table}

\begin{table}[H]
        \begin{center}
            \caption{Iteration numbers and CPU times of MINRES-$\mathcal{P}_{\alpha}$, MINRES-$\mathcal{P}_1$, and MINRES-$I$ with CN discretization for Example \ref{exampleheateq}.}\label{explcnhteqtbl_constant_heat}
            \setlength{\tabcolsep}{0.8em}
            {\small
                \begin{tabular}{ccc|cc|cc|cc}
                    \hline
                    \multirow{2}{*}{$N$}   & \multirow{2}{*}{$m+1$} & \multirow{2}{*}{DoF} & \multicolumn{2}{c|}{MINRES-$\mathcal{P}_{\alpha}$} & \multicolumn{2}{c|}{MINRES-$\mathcal{P}_{1}$} & \multicolumn{2}{c}{MINRES-$I$} \\ \cline{4-9} 
                    &                      &                      & Iter                    & CPU(s)                   & Iter                 & CPU(s)                 & Iter            & CPU(s)            \\ \hline
                    \multirow{4}{*}{$2^{5}$} & $2^{5}$             & 30752                & 2                       & 0.047                    & 66                   & 0.53                   & 35              & 0.092             \\
                    & $2^{6}$             & 127008               & 2                       & 0.056                     & 96                   & 1.38                   & 35              & 0.24              \\
                    & $2^{7}$             & 516128               & 2                       & 0.16                     & 141                  & 6.94                   & 35              & 0.70              \\
                    & $2^{8}$             & 2080800              & 2                       & 0.87                     & 184                  & 50.34                  & 36              & 3.23              \\ \hline
                    \multirow{4}{*}{$2^{6}$} & $2^{5}$             & 61504                & 2                       & 0.027                    & 95                   & 0.79                   & 64              & 0.33              \\
                    & $2^{6}$             & 254016               & 2                       & 0.086                    & 94                   & 2.25                   & 64              & 0.62              \\
                    & $2^{7}$             & 1032256              & 2                       & 0.35                     & 140                  & 13.78                  & 64              & 2.45              \\
                    & $2^{8}$             & 4161600              & 2                       & 2.16                     & 188                  & 108.90                 & 64              & 12.17             \\ \hline
                    \multirow{4}{*}{$2^{7}$} & $2^{5}$             & 123008               & 2                       & 0.050                    & 65                   & 1.01                   & 128             & 0.82              \\
                    & $2^{6}$             & 508032               & 2                       & 0.15                     & 94                   & 4.36                   & 128             & 2.19              \\
                    & $2^{7}$             & 2064512              & 2                       & 0.77                     & 139                  & 36.43                  & 128             & 10.06             \\
                    & $2^{8}$             & 8323200              & 2                       & 4.65                     & 226                  & 324.72                 & 128             & 76.47             \\ \hline
                    \multirow{4}{*}{$2^{8}$} & $2^{5}$             & 246016               & 2                       & 0.081                    & 65                   & 1.51                   & 256             & 2.34              \\
                    & $2^{6}$             & 1016064              & 2                       & 0.31                     & 93                   & 9.40                   & 256             & 10.07             \\
                    & $2^{7}$             & 4129024              & 2                       & 2.02                     & 139                  & 81.20                  & 256             & 48.74             \\
                    & $2^{8}$             & 16646400             & 2                       & 9.37                     & 228                  & 636.39                 & 256             & 297.46            \\ \hline
                \end{tabular}
            }
        \end{center}
\end{table}

}
\end{example}
										
In fact, MINRES-$\mathcal{P}_{\alpha}$ and MINRES-$\mathcal{P}_1$ are both applicable to the heat equation with variable coefficients \eqref{eqn:heat var coeff}. We demonstrate this via Example \ref{heateqvarcoeff}.
									
\begin{example}\label{heateqvarcoeff}
{\rm
In this example, we consider the heat equation \eqref{eqn:heat var coeff} with
\begin{align*}
    &d=2,~\Omega=(0,1)\times(0,1),~T=1,~u_0({ x})=x_2(1-x_2)x_1(1-x_1),\\
    &a({\bf x})=(20+ x_1^2)(20+ x_2^2),\\
    &f({ x},t)=\exp(t)[x_1(1-x_1)x_2(1-x_2)-2x_1(20+x_2^2)(1-2x_1)x_2(1-x_2)\\
    &\quad\quad\quad-2x_2(20+x_2^2)(1-2x_2)x_1(1-x_1)+2a({ x})(x_1(1-x_1)+x_2(1-x_2))].
\end{align*}
As BDF discretization is easier to handle, we here only consider the preconditioning techniques for CN discretization of Example \ref{heateqvarcoeff}. The preconditioning technique for BDF discretization is similar, which is thus neglected. The CN discretization of Example \ref{heateqvarcoeff} is as follows
\begin{equation}\label{eq:cnheatvarcoeff_tssystem}
  \frac{1}{\tau}  \begin{bmatrix}
        \check{L}_{a,1} &   &  & \\
        \check{L}_{a,2}  & 	\check{L}_{a,1}  & & \\
        & \ddots & \ddots &  \\
        &    & \check{L}_{a,2} & 	\check{L}_{a,1}
    \end{bmatrix}
    \left[\begin{array}{c}\mathbf{u}^{(1)}\\ \mathbf{u}^{(2)}\\ \vdots \\ \mathbf{u}^{(N)}
    \end{array} \right]
    = \left[\begin{array}{c}
        \mathbf{f}^{(\frac{1}{2})} -  \frac{1}{\tau}\check{L}_{a,2}{\bf u}^0 \\
        \mathbf{f}^{(\frac{3}{2})} \\
        \vdots \\
        \mathbf{f}^{(N-\frac{1}{2})}
    \end{array} \right],
\end{equation}
where $\check{L}_{a,1}=I-\frac{\tau}{2}\Delta_{a,h}$, $\check{L}_{a,2}=-I-\frac{\tau}{2}\Delta_{a,h}$; $\mathbf{f}^{(k-\frac{1}{2})}$ denotes the values of $f(\cdot,(k-1/2)\tau)$ on the spatial grid points for $k=1,2,...,N$. Clearly, the preconditioners $\mathcal{P}_{\alpha}$ and $\mathcal{P}_1$ cannot be applied to \eqref{eq:cnheatvarcoeff_tssystem} directly. However, if we consider replacing $a$ with its mean value $\bar{a}$ on the spatial grid, then we obtain another BLTT  matrix
\begin{equation}\label{meancoeffheateqtsmat}
   \bar{\mathcal{A}} = 
     \frac{1}{\tau}  \begin{bmatrix}
        
       \check{L}_{\bar{a},1} &   &  & \\
        \check{L}_{\bar{a},2}  & 	\check{L}_{\bar{a},1}  & & \\
        & \ddots & \ddots &  \\
        &    & \check{L}_{\bar{a},2} & 	\check{L}_{\bar{a},1}
    \end{bmatrix},
\end{equation}
where $\check{L}_{\bar{a},1}=I-\frac{\tau\bar{a}}{2}\Delta_{h}$, $\check{L}_{\bar{a},2}=-I-\frac{\tau\bar{a}}{2}\Delta_{h}$. The real symmetric blocks of the BLTT coefficient matrix $\mathcal{A}$ defined in \eqref{meancoeffheateqtsmat}  clearly satisfy Assumption \ref{akassumps},  due to the fact that $\Delta_h$ is RSND. Moreover, the preconditioners $\mathcal{P}_{\alpha}$ and $\mathcal{P}_1$ for \eqref{meancoeffheateqtsmat} are fast invertible, since $\check{L}_{\bar{a},1}$ and $\check{L}_{\bar{a},2}$ are diagonalizable by the sine transform. Hence, for Example \ref{heateqvarcoeff}, we consider the preconditioners $\mathcal{P}_{\alpha}$ and $\mathcal{P}_1$ defined for \eqref{meancoeffheateqtsmat} as preconditioners of preconditioned MINRES solvers for solving  the BLTT system \eqref{eq:cnheatvarcoeff_tssystem} (up to a symmetrization transformation). The numerical results of MINRES-$\mathcal{P}_{\alpha}$, MINRES-$\mathcal{P}_{1}$ and MINRES-$I$ for solving symmetrization of  \eqref{eq:cnheatvarcoeff_tssystem} are listed in Table \ref{explcnhteqvarcoeftbl_var_heat}. Table \ref{explcnhteqvarcoeftbl_var_heat} shows that MINRES-$\mathcal{P}_{\alpha}$ is the most efficient one among the three solvers in terms of both computation time and iteration number. 

\begin{table}[H]
        \begin{center}
            \caption{Iteration numbers and CPU times of MINRES-$\mathcal{P}_{\alpha}$, MINRES-$\mathcal{P}_1$, and MINRES-$I$ with CN discretization for Example \ref{heateqvarcoeff}.}\label{explcnhteqvarcoeftbl_var_heat}
            \setlength{\tabcolsep}{0.8em}
            {\small
                \begin{tabular}{ccc|cc|cc|cc}
                    \hline
                    \multirow{2}{*}{$N$}   & \multirow{2}{*}{$m+1$} & \multirow{2}{*}{DoF} & \multicolumn{2}{c|}{MINRES-$\mathcal{P}_{\alpha}$} & \multicolumn{2}{c|}{MINRES-$\mathcal{P}_{1}$} & \multicolumn{2}{c}{MINRES-$I$} \\ \cline{4-9} 
                    &                      &                      & Iter                    & CPU(s)                   & Iter                 & CPU(s)                 & Iter              & CPU(s)          \\ \hline
                    \multirow{4}{*}{$2^{5}$} & $2^{5}$             & 30752                & 10                      & 0.099                    & 217                  & 1.11                   & -             & -           \\
                    & $2^{6}$             & 127008               & 10                      & 0.18                     & 441                  & 5.68                   & -            & -          \\
                    & $2^{7}$             & 516128               & 10                      & 0.55                     & 893                  & 39.44                  & -           & -        \\
                    & $2^{8}$             & 2080800              & 10                      & 3.12                     & -                 & -                 & -         & -               \\ \hline
                    \multirow{4}{*}{$2^{6}$} & $2^{5}$             & 61504                & 10                      & 0.097                    & 184                  & 1.51                   & -             & -          \\
                    & $2^{6}$             & 254016               & 10                      & 0.29                     & 373                  & 8.47                   & -            & -         \\
                    & $2^{7}$             & 1032256              & 10                      & 1.20                     & 780                  & 71.73                  & -         & -               \\
                    & $2^{8}$             & 4161600              & 10                      & 7.43                     & -                 & -                 & -         & -               \\ \hline
                    \multirow{4}{*}{$2^{7}$} & $2^{5}$             & 123008               & 10                      & 0.15                     & 117                  & 1.51                   & -             & -          \\
                    & $2^{6}$             & 508032               & 10                      & 0.54                     & 234                  & 10.09                  & -         & -               \\
                    & $2^{7}$             & 2064512              & 10                      & 2.89                     & 471                  & 120.99                 & -         & -               \\
                    & $2^{8}$             & 8323200              & 10                      & 16.03                    & 949                  & 1424.86                & -         & -               \\ \hline
                    \multirow{4}{*}{$2^{8}$} & $2^{5}$             & 246016               & 10                      & 0.28                     & 65                   & 1.45                   & -             & -          \\
                    & $2^{6}$             & 1016064              & 10                      & 1.09                     & 129                  & 12.10                  & -         & -               \\
                    & $2^{7}$             & 4129024              & 10                      &  7.00                     & 257                  & 153.77                 & -         & -               \\
                    & $2^{8}$             & 16646400             & 10                      & 32.40                    & 516                  & 1511.59                & -         & -               \\ \hline
                \end{tabular}
            }
        \end{center}
\end{table}
}
\end{example}
\subsection{Non-local evolutionary equations}		To showcase the versatility of the proposed solver, we also utilize it to solve non-local evolutionary equations as presented below, even though this type of equation was not shown in the paper \cite{mcdonald2018preconditioning}.

Consider a non-local evolutionary equation with a weakly singular kernel \cite{lin2021parallel}:
\begin{equation}\label{eqn:evo eqn}
\left\{
\begin{aligned}
& \frac{1}{\Gamma(1-\gamma)} \int_0^{\rm T} \frac{\partial u(\mathbf{x}, s)}{\partial s}(t-s)^{-\gamma} d s=\nabla \cdot(a(\mathbf{x}) \nabla u)+f(\mathbf{x}, t), \quad \mathbf{x} \in \Omega \subset \mathbb{R}^2, t \in(0, T], \\
& u(\mathbf{x}, t)=0, \quad \mathbf{x} \in \partial \Omega, t \in(0, T], \\
& u(\mathbf{x}, 0)=u_0(\mathbf{x}), \quad \mathbf{x} \in \Omega,
\end{aligned}
\right.
\end{equation}
where $\Gamma(\cdot)$ is the gamma function, $\gamma \in(0,1), \Omega$ is open; $\partial \Omega$ denotes the boundary of $\Omega$;  $a$, $f$ and $u_0$ are all known functions.

With the L1 scheme \cite{lin2007finite,lin2018separable}, the temporal discretization has the following form
$$
\frac{1}{\Gamma(1-\gamma)} \int_0^{n \tau} \frac{\partial u(\mathbf{x}, s)}{\partial s}(t-s)^{-\gamma} d s \approx \frac{1}{\tau^\gamma} \sum_{k=1}^n l_{n-k}^{(\gamma)} u(\mathbf{x}, n \tau)+\frac{1}{\tau^\gamma} l^{(n, \gamma)} u_0(\mathbf{x}), \mathbf{x} \in \Omega, n=1,2, \ldots, N,
$$
where
$$
\begin{aligned}
& l_k^{(\gamma)}= \begin{cases}{[\Gamma(2-\gamma)]^{-1},} & k=0, \\
{[\Gamma(2-\gamma)]^{-1}\left[(k+1)^{1-\gamma}-2 k^{1-\gamma}+(k-1)^{1-\gamma}\right],} & 1 \leq k \leq N-1,\end{cases} \\
& l^{(k, \gamma)}=\left[(k-1)^{1-\gamma}-k^{1-\gamma}\right]\left[\tau^\gamma \Gamma(2-\gamma)\right]^{-1}, \quad 1 \leq k \leq N . \\
&
\end{aligned}
$$
We then obtain the discretization of \eqref{eqn:evo eqn} at each time step as follows
\begin{equation}\label{tsfulldisc}
\frac{1}{\tau^\gamma} \sum_{k=1}^n l_{n-k}^{(\gamma)} \mathbf{u}^k-\Delta_{a,h} \mathbf{u}^n=\mathbf{f}^n, \quad n=1,2, \ldots, N,
\end{equation}
where $\mathbf{u}^n \in \mathbb{R}^{M \times 1}$ is a vector whose components are approximate values of $u(\cdot, n \tau)$ on spatial grid points arranged in a lexicographic ordering, $\mathbf{f}^n$ contains the initial condition and the values of $f(\cdot, n \tau)$ on the spatial grid points.

Putting the $N$ many linear systems in an all-at-once linear system, we obtain the BLTT system as follows 
\begin{equation}\label{timefracblttsys}
\mathcal{A}{\bf u}={\bf f},
\end{equation}
where
\begin{align*}
&\mathbf{f}=\left(\mathbf{f}^{1}; \mathbf{f}^{2}; \cdots; \mathbf{f}^{N}\right)-\mathbf{v}_\gamma \otimes \mathbf{u}_0, \quad \mathbf{v}_\gamma=\left(l^{(1, \gamma)}, l^{(2, \gamma)}, \ldots, l^{(N, \gamma)}\right)^{\mathrm{T}}, \\
&\mathcal{A}=\mathbf{I}_N \otimes (-\Delta_{a,h})+\mathbf{T} \otimes \mathbf{I}_M, \quad \mathbf{T}=\frac{1}{\tau^\gamma}\left[\begin{array}{cccc}
	l_0^{(\gamma)} & & & \\
	l_1^{(\gamma)} & l_0^{(\gamma)} & & \\
	\vdots & \ddots & \ddots & \\
	l_{N-1}^{(\gamma)} & \ldots & l_1^{(\gamma)} & l_0^{(\gamma)}
\end{array}\right],\\
&\mathbf{u}=\left(\mathbf{u}^{1}; \mathbf{u}^2; \cdots; \mathbf{u}^N\right).
\end{align*}
$f^{n}$ denotes the values of $f(\cdot,n\tau)$ on the spatial grid points for $n=1,2,...,N$; the unknown vector ${\bf u}^{n}$ consists of approximation of the values of $u(\cdot,n\tau)$ on the spatial grid points for $n=1,2,...,N$, respectively.
Since ${\bf T}$ is a lower triangular Toeplitz matrix, $\mathcal{A}$ is actually a BLTT matrix of form \eqref{amatdef} with $A_{(0)}=\frac{l_0^{(\gamma)}}{\tau^\gamma}\mathbf{I}_M-\Delta_{a,h}$, $A_{(k)}=\frac{l_k^{(\gamma)}}{\tau^\gamma}\mathbf{I}_M $ for $k=1...N-1$. Also, the real symmetric blocks of the BLTT matrix $\mathcal{A}$ fulfill the Assumption \ref{akassumps}, since $-\Delta_{a,h}$ is RSND and the matrix ${\bf T}$ is a strictly diagonally dominant matrix with positive diagonal entries \cite{lin2021all,lin2018separable}.
\begin{example}\label{exam:3}
In this example, we consider the non-local evolutionary equation \eqref{eqn:evo eqn} with
$$
\begin{aligned}
& \Omega=(0, \pi) \times(0, \pi), \quad T=1, \quad \mathbf{x}=(x, y), \quad u(x, y, t)=\sin (x) \sin (y) t^2, \\
& a(x, y) \equiv 1, \quad f(x, y, t)=\sin (x) \sin (y)\left[\frac{2 t^{2-\gamma}}{\Gamma(3-\gamma)}+2 t^2\right].
\end{aligned}
$$	
\end{example}
Notice that Example \ref{exam:3} has a constant coefficient $a(\mathbf{x}) \equiv 1$. We test MINRES-$P_{\alpha}$, MINRES-$P_{1}$ and MINRES-$I$, the results of which are listed in Table \ref{tab:ex3r=0.1} ($\gamma=0.1$), Table \ref{tab:ex3r=0.5} ($\gamma=0.5$) and Table \ref{tab:ex3r=0.9} ($\gamma=0.9$), respectively. The numerical outcomes we have obtained lend support to the theoretical framework put forth earlier. Tables \ref{tab:ex3r=0.1} - \ref{tab:ex3r=0.9} show that (i) MINRES-$\mathcal{P}_{\alpha}$ is more efficient than MINRES-$\mathcal{P}_{1}$ in terms of computational time and iteration number; (ii) the convergence rate of MINRES-$\mathcal{P}_{\alpha}$ keeps stable as the temporal and the spatial grid refined. Notably, the number of iterations required for MINRES-$\mathcal{P}_{\alpha}$ remains stable as the temporal grid is refined, while the number of iterations required for MINRES-$\mathcal{P}_{1}$ increases.
Overall, the proposed MINRES-$\mathcal{P}_{\alpha}$ solver performs the best among the three solvers. 
\begin{table}[H]
		\begin{center}
			\caption{Iteration numbers and CPU times of MINRES-$\mathcal{P}_{\alpha}$, MINRES-$\mathcal{P}_1$, and MINRES-$I$ with $\gamma = 0.1$ for Example \ref{exam:3}.}\label{tab:ex3r=0.1}
			\setlength{\tabcolsep}{0.8em}
			{\small
			\begin{tabular}{ccc|cc|cc|cc}
			\hline
			\multirow{2}{*}{$N$}   & \multirow{2}{*}{$m+1$} & \multirow{2}{*}{DoF} & \multicolumn{2}{c|}{MINRES-$\mathcal{P}_{\alpha}$} & \multicolumn{2}{c|}{MINRES-$\mathcal{P}_{1}$} & \multicolumn{2}{c}{MINRES-$I$} \\ \cline{4-9} 
			&                      &                      & Iter                    & CPU(s)                   & Iter                 & CPU(s)                 & Iter            & CPU(s)            \\ \hline
			\multirow{4}{*}{$2^{5}$} & $2^{5}$             & 30752                & 2                       & 0.02                    & 6                   & 0.35                   & 180             & 0.44              \\
            & $2^{6}$             & 127008               & 2                       & 0.05                    & 6                   & 0.16                   & 696              & 4.66             \\
			& $2^{7}$             & 516128               & 2                       & 0.19                     & 6                  & 0.47                  &-              & -             \\
			& $2^{8}$             & 2080800              & 2                       & 1.07                     & 6                  & 2.25                  & -             & -             \\ \hline
		  \multirow{4}{*}{$2^{6}$} & $2^{5}$             & 61504                & 2                       & 0.03                    & 6                   & 0.08                   & 298              & 1.80              \\
			& $2^{6}$             & 254016               & 2                       & 0.09                    & 6                   & 0.22                  & -              & -              \\
			& $2^{7}$             & 1032256              & 2                       & 0.39                     & 6                  &  0.91                  & -              & -              \\
			& $2^{8}$             & 4161600              & 2                       & 2.47                     & 6                  & 5.14                 & -              & -             \\ \hline
			\multirow{4}{*}{$2^{7}$} & $2^{5}$             & 123008               & 2                       & 0.05                    & 6                   & 0.11                   & 298             &2.32              \\
			& $2^{6}$             & 508032               & 2                       & 0.17                     & 6                   & 0.38                   & -             & -              \\
			& $2^{7}$             & 2064512              & 2                       & 0.87                     & 6                  & 2.20                  & -             & -              \\
			& $2^{8}$             & 8323200              & 2                       & 5.04                     & 6                  & 10.68                 & -             & -             \\ \hline
			\multirow{4}{*}{$2^{8}$} & $2^{5}$             & 246016               & 2                       & 0.09                 & 8                   & 0.26                   & 322             & 3.35              \\
			& $2^{6}$             & 1016064              & 2                       & 0.37                     & 8                   & 1.10                   & -             & -              \\
			& $2^{7}$             & 4129024              & 2                       & 2.35                     & 8                  & 6.04                  & -             & -             \\
			& $2^{8}$             & 16646400             & 2                       & 9.75                     & 8                 & 25.73                 & -             & -            \\ \hline
			\end{tabular}
   }
		\end{center}
\end{table}
\begin{table}[H]
		\begin{center}
			\caption{Iteration numbers and CPU times of MINRES-$\mathcal{P}_{\alpha}$, MINRES-$\mathcal{P}_1$, and MINRES-$I$ with $\gamma = 0.5$ for Example \ref{exam:3}.}\label{tab:ex3r=0.5}
			\setlength{\tabcolsep}{0.8em}
			{\small
			\begin{tabular}{ccc|cc|cc|cc}
			\hline
			\multirow{2}{*}{$N$}   & \multirow{2}{*}{$m+1$} & \multirow{2}{*}{DoF} & \multicolumn{2}{c|}{MINRES-$\mathcal{P}_{\alpha}$} & \multicolumn{2}{c|}{MINRES-$\mathcal{P}_{1}$} & \multicolumn{2}{c}{MINRES-$I$} \\ \cline{4-9} 
			&                      &                      & Iter                    & CPU(s)                   & Iter                 & CPU(s)                 & Iter            & CPU(s)            \\ \hline
			\multirow{4}{*}{$2^{5}$} & $2^{5}$             & 30752                & 2                       & 0.02                    & 8                   & 0.07                   & 898              & 2.63              \\
            & $2^{6}$             & 127008               & 2                       & 0.05                    & 8                   & 0.15                  & -              & -              \\
			& $2^{7}$             & 516128               & 2                       & 0.19                     & 8                  & 0.50                   & -             & -              \\
			& $2^{8}$             & 2080800              & 2                       & 1.08                     & 8                  & 14.39                  & -              & -              \\ \hline
		  \multirow{4}{*}{$2^{6}$} & $2^{5}$             & 61504                & 2                       & 0.03                    & 8                   & 0.09                   & -              & -              \\
			& $2^{6}$             & 254016               & 2                       & 0.09                    & 8                   &  0.25                   & -              & -              \\
			& $2^{7}$             & 1032256              & 2                       & 0.39                     & 8                  & 1.03                  & -              & -              \\
			& $2^{8}$             & 4161600              & 2                       & 2.51                     & 8                  & 6.54                 & -              & -             \\ \hline
			\multirow{4}{*}{$2^{7}$} & $2^{5}$             & 123008               & 2                       & 0.05                    & 10                   & 0.16                   & -             & -              \\
			& $2^{6}$             & 508032               & 2                       & 0.18                     & 10                   & 0.57                   & -             & -              \\
			& $2^{7}$             & 2064512              & 2                       & 0.88                     & 10                  & 3.62                  & -             & -              \\
			& $2^{8}$             & 8323200              & 2                       & 5.06                     & 10                  & 16.53                 & -             & -             \\ \hline
			\multirow{4}{*}{$2^{8}$} & $2^{5}$             & 246016               & 2                       & 0.13                 & 10                   & 0.31                   & -             & -              \\
			& $2^{6}$             & 1016064              & 2                       & 0.36                     & 10                   & 1.14                   & -             & -              \\
			& $2^{7}$             & 4129024              & 2                       & 2.37                     & 10                  & 7.53                 & -             & -             \\
			& $2^{8}$             & 16646400             & 2                       & 9.77                     & 10                  & 31.32                 & -             & -            \\ \hline
			\end{tabular}
   }
		\end{center}
\end{table}


\begin{table}[H]
		\begin{center}
			\caption{Iteration numbers and CPU times of MINRES-$\mathcal{P}_{\alpha}$, MINRES-$\mathcal{P}_1$, and MINRES-$I$ with $\gamma = 0.9$ for Example \ref{exam:3}.}\label{tab:ex3r=0.9}
			\setlength{\tabcolsep}{0.8em}
			{\small
			\begin{tabular}{ccc|cc|cc|cc}
			\hline
			\multirow{2}{*}{$N$}   & \multirow{2}{*}{$m+1$} & \multirow{2}{*}{DoF} & \multicolumn{2}{c|}{MINRES-$\mathcal{P}_{\alpha}$} & \multicolumn{2}{c|}{MINRES-$\mathcal{P}_{1}$} & \multicolumn{2}{c}{MINRES-$I$} \\ \cline{4-9} 
			&                      &                      & Iter                    & CPU(s)                   & Iter                 & CPU(s)                 & Iter            & CPU(s)            \\ \hline
			\multirow{4}{*}{$2^{5}$} & $2^{5}$             & 30752                & 2                       & 0.02                    & 8                   &  0.06                   & 370              & 0.93              \\
            & $2^{6}$             & 127008               & 2                       & 0.05                   & 8                   & 0.14                   & -              & -             \\
			& $2^{7}$             & 516128               & 2                       & 0.20                     & 8                  & 0.50                   & -              & -              \\
			& $2^{8}$             & 2080800              & 2                       & 1.08                     & 8                  & 2.90                  & -              & -              \\ \hline
		  \multirow{4}{*}{$2^{6}$} & $2^{5}$             & 61504                & 2                       & 0.02                    & 8                   & 0.09                   & 568              & 3.41              \\
			& $2^{6}$             & 254016               & 2                       & 0.09                    & 8                   & 0.25                   & -              & -              \\
			& $2^{7}$             & 1032256              & 2                       & 0.40                     & 8                  & 1.04                  & -              & -              \\
			& $2^{8}$             & 4161600              & 2                       & 2.50                     & 8                  & 6.38                 & -              & -             \\ \hline
			\multirow{4}{*}{$2^{7}$} & $2^{5}$             & 123008               & 2                       & 0.05                   & 8                   &  0.14                   & 715             & 5.43              \\
			& $2^{6}$             & 508032               & 2                       & 0.17                     & 8                   & 0.46                   & -             & -              \\
			& $2^{7}$             & 2064512              & 2                       & 0.87                     & 8                  & 2.85                  & -             & -              \\
			& $2^{8}$             & 8323200              & 2                       & 5.07                     & 8                  & 13.47                 & -             & -             \\ \hline
			\multirow{4}{*}{$2^{8}$} & $2^{5}$             & 246016               & 2                       & 0.09                 & 10                   & 0.28                   & 884             & 9.32              \\
			& $2^{6}$             & 1016064              & 2                       & 0.36                     & 10                   & 1.13                   & -             & -              \\
			& $2^{7}$             & 4129024              & 2                       & 2.37                     & 10                  & 7.33                  & -             & -             \\
			& $2^{8}$             & 16646400             & 2                       & 9.73                     & 10                  & 31.11                 & -             & -            \\ \hline
			\end{tabular}
   }
		\end{center}
\end{table}

The rest of this section is devoted to testing efficiency of the proposed preconditioning technique on examples with non-constant $a(x, y)$.
\begin{example}\label{exam:4}
Consider the problem \eqref{eqn:evo eqn} with
\begin{equation*}
\begin{aligned}
& \Omega=(0,1)\times(0,1), \quad T=1, \quad \mathbf{x}=(x, y), \quad a(x, y)=35+x^{3.5}+y^{3.5}, \\
& u(x, y, t)=\sin (\pi x) \sin (\pi y) t^2, \quad f(x, y, t)=\sin (\pi x) \sin (\pi y)\left[\frac{2 t^{2-\alpha}}{\Gamma(3-\alpha)}+2 \pi^2 a t^2\right]\\
&-\pi t^2\left[\left(\partial_x a\right) \cos (\pi x) \sin (\pi y)+\left(\partial_y a\right) \sin (\pi x) \cos (\pi y)\right].
\end{aligned}
\end{equation*}	
\end{example}
 We test MINRES-$P_{\alpha}$, MINRES-$P_{1}$ and MINRES-$I$ on Example \ref{exam:4}, and the results of which are listed in Table \ref{tab:ex4r=0.3} ($\gamma=0.3$), Table \ref{tab:ex4r=0.1} ($\gamma=0.6$) and Table \ref{tab:ex4r=0.9} ($\gamma=0.9$), respectively. Tables \ref{tab:ex4r=0.1} - \ref{tab:ex4r=0.9} show that (i) MINRES-$\mathcal{P}_{\alpha}$ is more efficient than MINRES-$\mathcal{P}_{1}$ in terms of computational time and iteration number; (ii) the convergence rate of MINRES-$\mathcal{P}_{\alpha}$ does not deteriorate as the temporal and the spatial grid get refined. Additionally, we note that when $\gamma$ is relatively small, the iteration counts and CPU time for both methods are nearly identical. However, when $\gamma$ is large, the iteration counts of MINRES-$\mathcal{P}_{\alpha}$ and MINRES-$\mathcal{P}_{1}$ are similar, but the former outperforms the latter in terms of CPU time, requiring fewer resources to complete the computation. 

 \begin{table}[H]
		\begin{center}
			\caption{Iteration numbers and CPU times of MINRES-$\mathcal{P}_{\alpha}$, MINRES-$\mathcal{P}_1$, and MINRES-$I$ with $\gamma=0.3$ for Example \ref{exam:4}.}\label{tab:ex4r=0.3}
			\setlength{\tabcolsep}{0.8em}
			{\small
			\begin{tabular}{ccc|cc|cc|cc}
			\hline
			\multirow{2}{*}{$N$}   & \multirow{2}{*}{$m+1$} & \multirow{2}{*}{DoF} & \multicolumn{2}{c|}{MINRES-$\mathcal{P}_{\alpha}$} & \multicolumn{2}{c|}{MINRES-$\mathcal{P}_{1}$} & \multicolumn{2}{c}{MINRES-$I$} \\ \cline{4-9} 
			&                      &                      & Iter                    & CPU(s)                   & Iter                 & CPU(s)                 & Iter            & CPU(s)            \\ \hline
			\multirow{4}{*}{$2^{5}$} & $2^{5}$             & 30752                & 8                       &  0.05                    & 8                   & 0.06                   & -              & -              \\
            & $2^{6}$             & 127008               & 8                       & 0.14                    & 8                   & 0.13                   & -              & -              \\
			& $2^{7}$             & 516128               & 8                       & 0.45                     & 8                  & 0.45                   & -              & -              \\
			& $2^{8}$             & 2080800              & 8                       & 2.80                     & 8                  & 2.82                  & -              & -              \\ \hline
		  \multirow{4}{*}{$2^{6}$} & $2^{5}$             & 61504                & 8                       &  0.08                    & 8                   &  0.08                   & -              & -              \\
			& $2^{6}$             & 254016               & 8                       & 0.25                    & 8                   & 0.23                   & -              & -              \\
			& $2^{7}$             & 1032256              & 8                       &  1.05                     & 8                 &  1.05                  & -              & -              \\
			& $2^{8}$             & 4161600              & 8                       & 6.87                     & 8                  & 6.33                 & -              & -             \\ \hline
			\multirow{4}{*}{$2^{7}$} & $2^{5}$             & 123008               & 8                       &  0.12                    & 8                   & 0.13                   & -             & -              \\
			& $2^{6}$             & 508032               & 8                       & 0.47                     & 8                   & 0.45                   & -             & -              \\
			& $2^{7}$             & 2064512              & 8                       & 2.61                     & 8                  & 2.83                  & -             & -              \\
			& $2^{8}$             & 8323200              & 8                       &  13.79                     & 8                  & 14.36                 & -             & -             \\ \hline
			\multirow{4}{*}{$2^{8}$} & $2^{5}$             & 246016               & 8                       & 0.23                    & 8                   & 0.23                   & -             & -              \\
			& $2^{6}$             & 1016064              & 8                       & 0.94                     & 8                   & 1.12                   & -             & -              \\
			& $2^{7}$             & 4129024              & 8                       & 6.45                     & 8                  & 6.41                  & -             & -             \\
			& $2^{8}$             & 16646400             & 8                       & 26.98                     & 8                  & 26.95                 & -             & -            \\ \hline
			\end{tabular}
   }
		\end{center}
\end{table}

\begin{table}[H]
		\begin{center}
			\caption{Iteration numbers and CPU times of MINRES-$\mathcal{P}_{\alpha}$, MINRES-$\mathcal{P}_1$, and MINRES-$I$ with $\gamma=0.6$ for Example \ref{exam:4}.}\label{tab:ex4r=0.1}
			\setlength{\tabcolsep}{0.8em}
			{\small
			\begin{tabular}{ccc|cc|cc|cc}
			\hline
			\multirow{2}{*}{$N$}   & \multirow{2}{*}{$m+1$} & \multirow{2}{*}{DoF} & \multicolumn{2}{c|}{MINRES-$\mathcal{P}_{\alpha}$} & \multicolumn{2}{c|}{MINRES-$\mathcal{P}_{1}$} & \multicolumn{2}{c}{MINRES-$I$} \\ \cline{4-9} 
			&                      &                      & Iter                    & CPU(s)                   & Iter                 & CPU(s)                 & Iter            & CPU(s)            \\ \hline
			\multirow{4}{*}{$2^{5}$} & $2^{5}$             & 30752                & 8                       & 0.06                    & 8                   & 0.07                   & -              & -              \\
            & $2^{6}$             & 127008               & 8                       & 0.14                    & 8                   & 0.16                   & -              & -              \\
			& $2^{7}$             & 516128               & 8                       & 0.53                     & 8                  & 0.48                   & -              & -              \\
			& $2^{8}$             & 2080800              & 8                       & 2.92                     & 8                  & 2.90                  & -              & -              \\ \hline
		  \multirow{4}{*}{$2^{6}$} & $2^{5}$             & 61504                & 8                       & 0.08                    & 8                   &  0.09                   & -              & -              \\
			& $2^{6}$             & 254016               & 8                       & 0.26                    & 8                   & 0.24                   & -              & -              \\
			& $2^{7}$             & 1032256              & 8                       & 1.09                     & 8                  & 1.16                  & -              & -              \\
			& $2^{8}$             & 4161600              & 8                       & 6.89                     & 10                  & 7.89                 & -              & -             \\ \hline
			\multirow{4}{*}{$2^{7}$} & $2^{5}$             & 123008               & 8                       & 0.13                    & 8                   & 0.14                   & -             & -              \\
			& $2^{6}$             & 508032               & 8                       & 0.48                     & 10                   & 0.54                   & -             & -              \\
			& $2^{7}$             & 2064512              & 8                       & 2.61                    & 10                  & 3.42                 & -             & -              \\
			& $2^{8}$             & 8323200              & 8                       & 13.75                     & 10                  & 17.14                 & -             & -             \\ \hline
			\multirow{4}{*}{$2^{8}$} & $2^{5}$             & 246016               & 8                       & 0.24                    & 9                   & 0.24                   & -             & -              \\
			& $2^{6}$             & 1016064              & 8                       & 0.94                     & 10                   &  1.23                   & -             & -              \\
			& $2^{7}$             & 4129024              & 8                       & 6.47                     & 10                  & 7.50                  & -             & -             \\
			& $2^{8}$             & 16646400             & 8                       & 27.14                     & 10                  & 32.31                 & -             & -            \\ \hline
			\end{tabular}
   }
		\end{center}
\end{table}

\begin{table}[H]
		\begin{center}
			\caption{Iteration numbers and CPU times of MINRES-$\mathcal{P}_{\alpha}$, MINRES-$\mathcal{P}_1$, and MINRES-$I$ with $\gamma=0.9$ for Example \ref{exam:4}.}\label{tab:ex4r=0.9}
			\setlength{\tabcolsep}{0.8em}
			{\small
			\begin{tabular}{ccc|cc|cc|cc}
			\hline
			\multirow{2}{*}{$N$}   & \multirow{2}{*}{$m+1$} & \multirow{2}{*}{DoF} & \multicolumn{2}{c|}{MINRES-$\mathcal{P}_{\alpha}$} & \multicolumn{2}{c|}{MINRES-$\mathcal{P}_{1}$} & \multicolumn{2}{c}{MINRES-$I$} \\ \cline{4-9} 
			&                      &                      & Iter                    & CPU(s)                   & Iter                 & CPU(s)                 & Iter            & CPU(s)            \\ \hline
			\multirow{4}{*}{$2^{5}$} & $2^{5}$             & 30752                & 8                       & 0.06                    & 10                   &  0.08                   & -              & -              \\
            & $2^{6}$             & 127008               & 8                       & 0.15                    & 10                   & 0.16                   & -              & -              \\
			& $2^{7}$             & 516128               & 8                       & 0.53                     & 10                  & 0.59                   & -              & -              \\
			& $2^{8}$             & 2080800              & 8                       & 2.93                     & 10                  & 3.53                  & -              & -              \\ \hline
		  \multirow{4}{*}{$2^{6}$} & $2^{5}$             & 61504                & 8                       & 0.10                    & 10                   & 0.11                   & -              & -              \\
			& $2^{6}$             & 254016               & 8                       &  0.26                    & 10                   &  0.29                   & -              & -              \\
			& $2^{7}$             & 1032256              & 8                       & 1.11                     & 10                  & 1.4                  & -              & -              \\
			& $2^{8}$             & 4161600              & 8                       & 6.90                     & 10                  & 7.87                 & -              & -             \\ \hline
			\multirow{4}{*}{$2^{7}$} & $2^{5}$             & 123008               & 8                       & 0.13                    & 10                   & 0.17                    & -             & -              \\
			& $2^{6}$             & 508032               & 8                       & 0.49                     & 10                   & 0.54                   & -             & -              \\
			& $2^{7}$             & 2064512              & 8                       & 2.62                     & 10                  & 3.41                  & -             & -              \\
			& $2^{8}$             & 8323200              & 8                       & 13.91                     & 10                  &  17.22                 & -             & -             \\ \hline
			\multirow{4}{*}{$2^{8}$} & $2^{5}$             & 246016               & 8                       & 0.24                    & 10                   & 0.28                   & -             & -              \\
			& $2^{6}$             & 1016064              & 8                       & 0.94                     & 10                   & 1.29                   & -             & -              \\
			& $2^{7}$             & 4129024              & 8                       & 6.49                      & 10                  & 7.5                  & -             & -             \\
			& $2^{8}$             & 16646400             & 8                       & 26.97                     & 10                  & 32.28                 & -             & -            \\ \hline
			\end{tabular}
   }
		\end{center}
\end{table}

\section{Conclusions}\label{conclusions}

In this paper, we have generalized the absolute value block circulant preconditioner proposed in \cite{mcdonald2018preconditioning} for the symmetrized BLTT system \eqref{symmetrizedblttsys} to the ABAC preconditioner by introducing a parameter $\alpha\in(0,1]$.  Fast implementation for the ABAC preconditioner has been discussed, which leads to a linearithmic complexity (nearly optimal) for each preconditioned MINRES iteration. Most importantly, with our proposed ABAC preconditioner, we have shown that under properly chosen  $\alpha$, the preconditioned MINRES solver has a matrix-size independent convergence rate for the symmetrized BLTT system. Numerical results reported have demonstrated the efficiency, versatility of the proposed preconditioning method and supported the theoretical results. 
\begin{acknowledgements}
The work of Sean Hon was supported in part by the Hong Kong RGC under grant 22300921, a start-up grant from the Croucher Foundation, and a Tier 2 Start-up Grant from Hong Kong Baptist University. The work of Xuelei Lin was supported by research grants: 2021M702281 from China Postdoctoral Science Foundation and HA45001143, a start-up Grant from Harbin Institute of Technology, Shenzhen. The work of Shu-Lin Wu was supported by Jilin Provincial Department of Science and Technology (No. YDZJ202201ZYTS593).
\end{acknowledgements}

\section*{Data Availability Statement}
The data of natural images and codes involved in this paper are available from the corresponding author on reasonable request.

\section*{Declarations}
\begin{description}
	\item[1.]I confirm that I have read, understand and agreed to the submission guidelines, policies and submission 
	declaration of the journal.\\
	\item[2.]I confirm that all authors of the manuscript have no conflict of interests to declare.\\
	\item[3.] I confirm that the manuscript is the authors' original work and the manuscript has not received prior publication 
	and is not under consideration for publication elsewhere.\\
	\item[4.] On behalf of all Co-Authors, I shall bear full responsibility for the submission.\\
	\item[5.] I confirm that all authors listed on the title page have contributed significantly to the work, have read the 
	manuscript, attest to the validity and legitimacy of the data and its interpretation, and agree to its submission.\\
	\item[6.]  I confirm that the paper now submitted is not copied or plagiarized version of some other published work.\\
	\item[7.] I declare that I shall not submit the paper for publication in any other Journal or Magazine till the decision is 
	made by journal editors. \\
	\item[8.] If the paper is finally accepted by the journal for publication, I confirm that I will either publish the paper
	immediately or withdraw it according to withdrawal policies.\\
	\item[9.] I understand that submission of false or incorrect information/undertaking would invite appropriate penal 
	actions as per norms/rules of the journal.
\end{description}

%
%

\bibliographystyle{spmpsci}      

\end{document}